\documentclass{article}
\usepackage{amssymb}
\usepackage{amsfonts}
\usepackage{graphicx}
\usepackage{amsmath}

\newtheorem{theorem}{Theorem}

\newtheorem{conclusion}[theorem]{Conclusion}

\newtheorem{example}[theorem]{Example}

\newtheorem{lemma}[theorem]{Lemma}

\newtheorem{problem}[theorem]{Problem}
\newtheorem{proposition}[theorem]{Proposition}
\newtheorem{remark}[theorem]{Remark}

\newenvironment{proof}[1][Proof]{\noindent\textbf{#1.} }{\ \rule{0.5em}{0.5em}}

\begin{document}

\title{\textsf{Topology of partition of measures by fans and the second
obstruction}}
\author{Pavle V.M. Blagojevi\'{c}\thanks{{\small Supported by St. John the
Baptist Serbian Orthodox Church in San Francisco, the grant 1854 of the
Serbian Ministry of Science, Technology and Development, by Town goverments
of Leskovac and Sremska Mitrovica and Water supply company of Leskovac.}} \\
Mathematical Institute SANU, Belgrade}
\maketitle

\begin{abstract}
\noindent The simultaneous partition problems are classical problems of the
combinatorial geometry which have the natural flavor of the equivariant
topology. The $k$-fan partition problems have attracted a lot of attention
\cite{Aki2000}, \cite{BaMa2001}, \cite{BaMa2002} and forced some hard
concrete combinatorial calculations in the equivariant cohomology \cite%
{Bl-Vr-Ziv}. These problems can be reduced, by a beautiful scheme of \cite%
{BaMa2001}, to a \textquotedblright typical\textquotedblright\ question of
the existence of a $\mathbb{D}_{2n}$ equivariant map $f:V_{2}(\mathbb{R}%
^{3})\rightarrow W_{n}-\cup \mathcal{A}(\alpha )$, \ where $V_{2}(\mathbb{R}%
^{3})$ is the space of all orthonormal $2$-frames in $\mathbb{R}^{3}$ and $%
W_{n}-\cup \mathcal{A}(\alpha )$ is the complement of the appropriate
arrangement. We introduce the \textit{target extension scheme} which allow
us to use the equivariant obstruction theory as a tool for proving that: for
every two proper measures on the sphere $S^{2}$, and any $\alpha
=(a,a+b,b)\in \mathbb{R}_{>0}^{3}$, there exists an $\alpha $-partition of
theses measures by a $3$-fan.

\noindent The significance of these results, among other, is that, beside
negative results \cite{Bl-Vr-Ziv}, the equivariant obstruction theory can
pull off some positive results, which were not attained by other means.
\end{abstract}

\section{\textsf{Introduction}}

\subsection{Problem}

\noindent A $k$-fan $(x;l_{1},l_{2},\ldots ,l_{k})$ on the sphere $S^{2}$ is
formed of a point $x$, called the center of the fan, and $k$ great
semicircles $l_{1},\ldots ,l_{k}$ emanating from $x$. We always assume
counter clockwise enumeration on great semicircles $l_{1},\ldots ,l_{k}$ of
a $k$-fan. Sometimes instead of lines we use open angular sectors $\sigma
_{i}$ between $l_{i}$ and $l_{i+1},\,i=1,\ldots ,k$. In that case we denote
a $k$-fan with $(x;\sigma _{1},\sigma _{2},\ldots ,\sigma _{k})$.

\noindent Let $\mu _{1},\mu _{2},\ldots ,\mu _{m}$ be \textit{proper} Borel
probability measures on $S^{2}$. Measure $\mu $ is \textit{proper} if $\mu
([a,b])=0$ for any circular arc $[a,b]\subset S^{2}$ and $\mu (U)>0$ for
each nonempty open set $U\subset S^{2}$. All results can be extended to more
general measures, including the counting measures of finite sets, see \cite%
{BaMa2001} for related examples.

\noindent Let $(\alpha _{1},\alpha _{2},\ldots ,\alpha _{k})\in \mathbb{R}%
_{>0}^{k}$ be a vector where $\alpha _{1}+\alpha _{2}+\ldots +\alpha _{k}=1$%
. The general problem stated in \cite{BaMa2001} is:

\begin{problem}
Find all triples $(m,k,\alpha )\in \mathbb{N\times }\mathbb{N\times R}^{k}$
such that for any collection of $m$ measures $\{\mu _{1},\mu _{2},\ldots
,\mu _{m}\}$, there exists a $k$-fan $(x;l_{1},\ldots ,l_{k})$ with the
property%
\begin{equation*}
(\forall i=1,\ldots ,k)\,(\forall j=1,\ldots ,m)\,\mu _{j}(\sigma
_{i})=\alpha _{i}
\end{equation*}%
That kind of a $k$-fan $(x;l_{1},\ldots ,l_{k})$ is called an $\alpha $%
\textit{-partition} for the collection of measures $\{\mu _{j}\}_{j=1}^{m}$.
\end{problem}

\noindent The analysis given in \cite{BaMa2001} shows that the most
interesting triples are $(3,2,\alpha ),(2,3,\alpha ),(2,4,\alpha )$. Known
results can be summed in the following table:%
\begin{equation*}
\begin{tabular}{|l|l|l|l|}
\hline
$m/k$ & $2$ & $3$ & $4$ \\ \hline
$2$ & $R_{\geq 0}^{2},\text{\cite{BaMa2001}}$ & $%
\begin{array}{c}
(\tfrac{1}{3},\tfrac{1}{3},\frac{1}{3}) \\
(\frac{1}{2},\frac{1}{4},\frac{1}{4}) \\
(\frac{\ast }{5},\frac{\ast }{5},\frac{\ast }{5})%
\end{array}%
\text{\cite{BaMa2001}}$ & $%
\begin{array}{c}
(\frac{2}{5},\frac{1}{5},\frac{1}{5},\frac{1}{5})\text{\cite{BaMa2001}} \\
(\frac{1}{4},\frac{1}{4},\frac{1}{4},\frac{1}{4})\text{\cite{BaMa2002}} \\
\text{{\small Topology can't help any more \cite{Bl-Vr-Ziv}.}}%
\end{array}%
$ \\ \hline
$3$ & $%
\begin{array}{c}
(\frac{1}{2},\frac{1}{2}) \\
(\frac{1}{3},\frac{2}{3}) \\
(\frac{2}{3},\frac{1}{3})%
\end{array}%
\text{\cite{BaMa2001}}$ & $\emptyset $ & $\emptyset $ \\ \hline
\end{tabular}%
\
\end{equation*}%
In this paper we try to fill this table a little bit more.

\subsection{The statement of results}

\noindent We are interested in the problem of $3$-fan partitions of two
measures on $S^{2}$, respectively. We prove the following result:

\begin{theorem}
\label{th:main1}Lets choose $\alpha =(a,a+b,b)\in \mathbb{R}_{>0}^{3}$ such
that $2a+2b=1$. Then any two \textit{proper} measures $\mu $ and $\nu $ on
the sphere $S^{2}$ admit an $\alpha $-partition by a $3$-fan $\mathfrak{p}%
=(x;l_{1},l_{2},l_{3})$, i. e.
\begin{equation*}
\mathcal{A}_{2,3}\supseteq \{(a,a+b,b)\in \mathbb{R}^{3}~|~a,b>0,~2a+b=1\}.
\end{equation*}
\end{theorem}

\begin{remark}
The case $a=b=1$ was already considered in \cite{BaMa2001}
\end{remark}

\subsection{The solution scheme}

\noindent The solution of the problem has two natural parts. The first part
is the reduction of the problem to the question of the existence of the
appropriate equivariant map. The second part is and topological effort to
disprove the existence of a such map. There is also a third part of the
proof, the limit argument. It extends the result from rational triples to
real triples, but we omited it because it is the standard part of every
similiar proof.

\medskip

\noindent \textit{Reduction to the equivariant problem:}

\begin{itemize}
\item The configuration space / test map procedure of Imre B\'{a}r\'{a}ny
and Ji\v{r}i Matou\v{s}ek from \cite{BaMa2001} reduces the problem to the
question: Is there an $\alpha =(\frac{a_{1}}{n}$,$\frac{a_{2}}{n}$,$\frac{%
a_{3}}{n})$, such that there is no $\mathbb{D}_{2n}$-map $V_{2}(\mathbb{R}%
^{3})\rightarrow W_{n}\setminus \cup \mathcal{A}(\alpha )$?

\item The extension of scalars\ equivalence from homological algebra, \cite%
{Brown}, allows us to change the initial equivariant question to: Is there
an $\alpha =(\frac{a_{1}}{n},\frac{a_{2}}{n},\frac{a_{3}}{n})$, such that
there is no $\mathbb{Q}_{4n}$-map $S^{3}\rightarrow W_{n}\setminus \cup
\mathcal{A}(\alpha )$?

\item The elementary obstruction theory can convince us that the actual free
$\mathbb{Q}_{4n}$ action on $S^{3}$ is not of the essential importance. Thus
we can always change it the way it pleases us.
\end{itemize}

\medskip

\noindent \textit{Obstruction theory approach:}

\begin{itemize}
\item The dimensional reasons of the problem, $\pi _{1}(W_{n}\setminus \cup
\mathcal{A}(\alpha ))\neq \varnothing $, implies that there are two
essential obstructions. This forces us to introduce the target extension
scheme, and instead of dealing with two obstructions we have only one to
compute. The target extension scheme extends the $\mathbb{Q}_{4n}$ space $%
W_{n}\setminus \cup \mathcal{A}(\alpha )$ in such a way that connectivity
increases by one. Thus, the problem of the existence of  the $\mathbb{Q}_{4n}
$-map from $S^{3}$ to this extension has only one obstruction.

\item Once we extend the complement $W_{n}\setminus \cup \mathcal{A}(\alpha )
$ we use the classical obstruction theory and the "map in general position"
method to prove that appropriate obstruction is not zero. This part of the
proof goes in a small number of steps, like we described in the section \ref%
{sec:HowToComputeTheOstructionCocycle}.

\item The main step is the identification of the obstruction cocycle. This
is the point where many previous papers have been broken (\cite{VreZiv2002},
proof of theorem 4.2, equality (25) and \cite{limun2002}, proof of theorem
6.1). We use the clear geometrical picture and simple testing methods (in
sections \ref{Sec:Step7} and \ref{sec:Step8}) to interpret the obstruction
element as a non-zero element of the appropriate coinvariant group.
\end{itemize}

\noindent This scheme of the solution, with different extensions, can be
applied to other cases of the $3$-fan / $2$-measures problem, as well on the
$2$-fan / $3$-measures problem. The main idea of extending the target space,
of course with some variations, can be applied on every problem of the
existace of the equivaiant map to a complement of an arrangement, or to any
space that has natural extension candidates.

\section{\textsf{From a partition problem to an equivariant problem}}

\subsection{Reduction to the equivariant problem}

\noindent The reduction of the fan partition problem to the equivariant
problem is done by the configurations space / test map scheme. The main idea
is to look at the space of all possible solutions and to rephrase the
question in terms of coincidences of the associated test map. Imre B\'{a}r%
\'{a}ny and Ji\v{r}i Matou\v{s}ek demonstrated in \cite{BaMa2001} that the
test map scheme can be applied on the problem of $\alpha $-partitioning of $%
m $-measures on $S^{2}$ by spherical $k$-fans. In a very elegant way this
problem was reduced to the problem of the existence of the appropriate
equivariant map. We briefly review this reduction for the $(3,2)$ case of
this problem.

\noindent \textbf{The configuration space.} Let $\mu $ and $\nu $ be two
proper Borel probability measures on $S^{2}$, and $F_{k}$ the space of all $%
k $-fans on the sphere $S^{2}$. The space $X_{\mu }$ of all possible
solutions associated to the measure $\mu $ is defined by
\begin{equation*}
X_{\mu }=\{(x;l_{1},\ldots ,l_{n})\in F_{n}\mid (\forall i=1,\ldots ,n)\,\mu
(\sigma _{i})=\tfrac{{\small 1}}{{\small n}}\}.
\end{equation*}%
Observe that every $n$-fan $(x;l_{1},\ldots ,l_{n})\in X_{\mu }$ is
completely determined by the pair $(x,l_{1})$ or equivalently, the pair $%
(x,y)$, where $y$ is the unit tangent vector to $l_{1}$ at $x$. Thus, the
space $X_{\mu }$ is a Stiefel manifold $V_{2}(\mathbb{R}^{3})$ of all
orthonormal $2$-frames in $\mathbb{R}^{3}$. Keep in mind that $V_{2}(\mathbb{%
R}^{3})\cong SO(3)\cong \mathbb{R}P^{3}$.

\noindent \textbf{Test map.} Let $\mathbb{R}^{n}$ be an Euclidean
space with the standard orthonormal basis $e_{1},e_{2},\ldots
,e_{n}$ and the associated coordinate functions
$x_{1},x_{2},\ldots ,x_{n}$. Let $W_{n}$ be
the hyperplane $\{x\in \mathbb{R}^{n}\mid x_{1}+x_{2}+\ldots +x_{n}=0\}$ in $%
\mathbb{R}^{n}$, and suppose that $\alpha $-vectors have the following form
\begin{equation*}
\alpha =(\tfrac{a_{1}}{n},\tfrac{a_{2}}{n},\tfrac{a_{3}}{n})\in \frac{1}{n}\,%
\mathbb{N}^{3}\subseteq \mathbb{Q}^{3},
\end{equation*}%
where $a_{1}+a_{2}+a_{3}=n$. Then test maps for the $(3,2)$ fan problem are
defined by%
\begin{equation*}
\begin{array}{lll}
F_{\nu }:X_{\mu }\rightarrow W_{n} &  & F_{\nu }(\mathfrak{p})=(\nu (\sigma
_{1})-\tfrac{{\small 1}}{{\small n}},\ldots ,\nu (\sigma _{n})-\tfrac{%
{\small 1}}{{\small n}})%
\end{array}%
\end{equation*}%

\noindent \textbf{The action.} The dihedral group
$\mathbb{D}_{2n}=\langle j,\varepsilon \,|\,\varepsilon
^{n}=j^{2}=1,\,\varepsilon j=j\varepsilon ^{n-1}\,\rangle $ acts
both on the possible solution space $X_{\mu }$ and the linear
subspace $W_{n}\subseteq \mathbb{R}^{n}$ by
\begin{equation*}
X_{\mu }:\left\{
\begin{array}{l}
\varepsilon (x;l_{1},\ldots ,l_{n})=(x;l_{n},l_{1},\ldots ,l_{n-1}) \\
j(x;l_{1},\ldots ,l_{n})=(-x;l_{1},l_{n},l_{n-1,}\ldots ,l_{2})%
\end{array}%
,\right. W_{n}:\left\{
\begin{array}{c}
\varepsilon (x_{1},\ldots ,x_{n})=(x_{2},\ldots ,x_{n},x_{1}) \\
j(x_{1},\ldots ,x_{n})=(x_{n},\ldots ,x_{2},x_{1})%
\end{array}%
,\right.
\end{equation*}%
for $(x;l_{1},\ldots ,l_{n})\in X_{\mu }$ and $(x_{1},\ldots
,x_{n})\in W_{n} $. The action of $\mathbb{D}_{2n}$ on $X_{\mu }$
is \textbf{free}.

\noindent Observe that the space of possible solutions $X_{\mu }$ is $%
\mathbb{D}_{2n}$-homeomorphic to the manifold
$V_{2}(\mathbb{R}^{3})$, where $V_{2}(\mathbb{R}^{3})$ is a
$\mathbb{D}_{2n}$-space given by
\begin{equation*}
\varepsilon (x,y)=(x,R_{x}({\frac{2\pi }{n}})(y))\text{, }j(x,y)=(-x,y)\text{%
,}
\end{equation*}%
and $R_{x}(\theta ):\mathbb{R}^{3}\rightarrow \mathbb{R}^{3}$ is
the rotation round the axes determined by $x$ through the angle
$\theta $.

\noindent \textbf{The test space. }The test space in this problem
is the union $\cup \mathcal{A}(\alpha )\subset W_{n}$ of a
smallest $\mathbb{D}_{2n} $-invariant linear subspace arrangement
$\mathcal{A}(\alpha )$, which contains linear subspace $L(\alpha
)\subset W_{n}$. The subspace $L(\alpha )$ is defined by
\begin{equation*}
L(\alpha )=\{x\in \mathbb{R}^{n}\mid \xi _{1}(x)=\xi _{2}(x)=\xi
_{3}(x)=0\}\subseteq W_{n},
\end{equation*}%
where
\begin{equation*}
\begin{array}{lll}
\xi _{1}(x)=x_{1}+\ldots +x_{a_{1}}, & \xi
_{2}(x)=x_{a_{1}+1}+\ldots
+x_{a_{1}+a_{2}}, & \xi _{3}(x)=x_{a_{1}+a_{2}+1}+\ldots +x_{n}\text{.}%
\end{array}%
\end{equation*}

\noindent Since the test map $F_{\nu }$ is obviously $\mathbb{D}_{2n}$%
-equivariant, the following standard proposition is proved.

\begin{proposition}
\label{prop:VezaProblem-Ekvivarijantan}Let $\alpha =(\frac{a_{1}}{n},\frac{%
a_{2}}{n},\frac{a_{3}}{n})\in \frac{1}{n}\,\mathbb{N}^{3}\subseteq \mathbb{Q}%
^{3}$ be a vector such that $a_{1}+a_{2}+a_{3}=n$. If there is no $\mathbb{D}%
_{2n}$-equivariant map%
\begin{equation*}
F:V_{2}(\mathbb{R}^{3})\rightarrow W_{n}\setminus \cup \mathcal{A}(\alpha )
\end{equation*}%
then for any two measures $\mu $ and $\nu $ on $S^{2}$, there exists an $%
\alpha $-partition $(x;l_{1},l_{2},l_{3})$ of measures $\mu $ and $\nu $.
\end{proposition}

\subsection{Modifying Problem}

\noindent Knowing the fact that for an odd $n$, there always exists a $%
\mathbb{Z}_{n}$-map $f:S^{3}\rightarrow V_{2}(\mathbb{R}^{3})$, B\'{a}r\'{a}%
ny and Matou\v{s}ek in \cite{BaMa2001} questioned if there is a $\mathbb{Z}%
_{n}$-map from the sphere $S^{3}$ to the complements of appropriate
arrangements. To do something similar we extend the group, like in \cite%
{limun2002} and \cite{Bl-Vr-Ziv}. We use well known \textquotedblleft
extension of scalars\textquotedblright\ equivalence from homological
algebra, \cite{Brown} Section III.3.

\noindent \textbf{The generalized quaternion group.} Let $S^{3}=S(\mathbb{H}%
)=Sp(1)$ be the group of all unit quaternions and let $\epsilon =\epsilon
_{2n}=\cos \frac{\pi }{n}+i\sin \frac{\pi }{n}\in S(\mathbb{H})$ be a root
of unity. Group $\langle \epsilon \rangle $ is a subgroup of $S(\mathbb{H})$
of the order $2n$. Then, the generalized quaternion group\textit{, }\cite%
{CaEi}\textit{\ }p. 253\textit{, }is the subgroup
\begin{equation*}
\mathbb{Q}_{4n}=\{1,\epsilon ,\ldots ,\epsilon ^{2n-1},j,\epsilon j,\ldots
,\epsilon ^{2n-1}j\}
\end{equation*}%
of $S^{3}$, of the order $4n$. Let $H=\{1,\epsilon ^{n}\}=\{1,-1\}\subset
\mathbb{Q}_{4n}$. Then, it is not hard to see that the quotient group $%
\mathbb{Q}_{4n}/H$ is isomorphic to the dihedral group $\mathbb{D}_{2n}$ of
the order $2n$.\textbf{\ }

\begin{proposition}
\label{prop:PrelazakNaNovuGrupu}Let the generalized quaternion group\textit{%
\ }$\mathbb{Q}_{4n}$ act on $S^{3}$ as a subgroup, and on $W_{n}$ via
already defined $\mathbb{D}_{2n}$ action by the quotient homomorphism $%
\mathbb{Q}_{4n}\rightarrow \mathbb{Q}_{4n}/\{1,-1\}\cong \mathbb{D}_{2n}$.
Then the following maps coexist:%
\begin{equation*}
\mathbb{D}_{2n}\text{-map }V_{2}(\mathbb{R}^{3})\rightarrow W_{n}\setminus
\cup \mathcal{A}(\alpha )\text{ \ and \ }\mathbb{Q}_{4n}\text{-map }%
S^{3}\rightarrow W_{n}\setminus \cup \mathcal{A}(\alpha ).
\end{equation*}

\noindent By the coexistence we mean that the one map exists if and only if
the other map exists, i.e. the one can't exist without the other.
\end{proposition}

\begin{proof}
Let us denote the target space $W_{n}\setminus \cup \mathcal{A}(\alpha )$
with $T$. Also, observe that $S^{3}/\{1,-1\}\cong \mathbb{R}P^{3}\cong
SO(3)\cong V_{2}(\mathbb{R}^{3})$

\noindent $\Rightarrow :$ Let $F:V_{2}(\mathbb{R}^{3})\rightarrow T$ be a $%
\mathbb{D}_{2n}$-map. The quotient map $p:S^{3}\rightarrow
S^{3}/\{1,-1\}\cong V_{2}(\mathbb{R}^{3})$ is a $\mathbb{Q}_{4n}$-map where
the $\mathbb{Q}_{4n}$-acts on $V_{2}(\mathbb{R}^{3})$ by the quotient
homomorphism $\mathbb{Q}_{4n}\rightarrow \mathbb{Q}_{4n}/\{1,-1\}\cong
\mathbb{D}_{2n}$. Since, $\mathbb{Q}_{4n}$ acts on both $V_{2}(\mathbb{R}%
^{3})$ and $T$ via the quotient homomorphism, the given $\mathbb{D}_{2n}$%
-map $F:V_{2}(\mathbb{R}^{3})\rightarrow T$ can also be seen as the $\mathbb{%
Q}_{4n}$-map. Thus, the composition
\begin{equation*}
F\circ p:S^{3}\longrightarrow S^{3}/\{1,-1\}\cong V_{2}(\mathbb{R}%
^{3})\longrightarrow T
\end{equation*}%
is the required $\mathbb{Q}_{4n}$-map $S^{3}\rightarrow T$.

\noindent $\Rightarrow :$ Let $G:S^{3}\rightarrow T$ be a $\mathbb{Q}_{4n}$%
-map. Observe that $S^{3}/\{1,-1\}\cong V_{2}(\mathbb{R}^{3})$ can be seen
as a $\mathbb{D}_{2n}\cong \mathbb{Q}_{4n}/\{1,-1\}$ space by $\{1,-1\}x%
\overset{g\{1,-1\}}{\longmapsto }\{1,-1\}(gx)$, where $x\in S^{3}$ and $g\in
\mathbb{Q}_{4n}$. Since the subgroup $\{1,-1\}$ acts trivially on $T$, there
is a factorization of the map $G$ through the quotient $h:S^{3}/\{1,-1\}%
\rightarrow T$ such that $G=h\circ p$. The map $h$ is the required $\mathbb{D%
}_{2n}$-map%
\begin{equation*}
h(g\{1,-1\}\cdot \{1,-1\}x)=G(g\cdot x)=g\cdot G(x)=g\cdot
h(\{1,-1\}x)=(g\{1,-1\})\cdot h(\{1,-1\}x).
\end{equation*}
\end{proof}

\begin{remark}
The $\mathbb{Q}_{4n}$ action on $S^{3}$ is \textbf{free}. Also, the $\mathbb{%
Q}_{4n}$ action on $W_{n}$ is the restriction of the following $\mathbb{Q}%
_{4n}$ action on $\mathbb{R}^{n}$. Let $e_{1},..,e_{n}$ be the standard
orthonormal basis in $\mathbb{R}^{n}$. The action is defined by
\begin{equation*}
\epsilon \cdot e_{i}=e_{i \ \textrm{mod}\ n+1} \text{ and }j\cdot e_{i}=e_{n-i+1}%
\text{.}
\end{equation*}
\end{remark}

\noindent \textbf{The free action on }$S^{3}$\textbf{.} Since the sphere $%
S^{3}$ is $2$-connected it turns out that the particular $\mathbb{Q}_{4n}$%
-action on $S^{3}$ is not something we have to live with. The elementary
equivariant obstruction theory allows us to prove the following useful fact.

\begin{proposition}
\label{prop:promenaDejstva}If $\gamma _{1}$ and $\gamma _{2}$ are $G$%
-actions on $S^{3}$ and $\gamma _{1}$ is free, then there exists a $G$-map $%
f:S^{3}\rightarrow S^{3}$ such that
\begin{equation*}
(\forall g\in G)\,(\forall x\in S^{3})\,\,f(g\cdot _{1}x)=g\cdot _{2}f(x)%
\text{.}
\end{equation*}
\end{proposition}

\begin{proof}
The statement is true because $S^{3}$ is $2$-connected, $\gamma _{1}$ is
free and so there are no obstructions to extend a $G$-map from $0$-skeleton
to $S^{3}$.
\end{proof}

\noindent Thus, when the time comes we will be able to use any free $\mathbb{%
Q}_{4n}$ action on $S^{3}$ and we will have the favorite one.

\section{\textsf{Obstruction theory approach}}

\noindent Once again, denote the target space $W_{n}\setminus \cup \mathcal{A%
}(\alpha )$ with $T$. To answer a question of the existence of $\mathbb{Q}%
_{4n}$-map $S^{3}\rightarrow T$ we will try to employ the classical
obstruction theory. Since the maximal elements of the arrangement $\mathcal{A%
}(\alpha )$ are of the codimension $2=(n-1)-(n-3)$ in $W_{n}$, the
complement is connected and the first obstruction lives in $H_{\mathbb{Q}%
_{4n}}^{2}(S^{3},\pi _{1}(T))$. In this case it is very hard even to
identify this group, not to mention to identify the particular element in
it. Even if we do manage to identify and calculate the first obstruction,
there is a good chance that it is zero, so the second obstructions should be
calculated. To omit this difficulties we use the nature of the target
spaces, introduce the target extension scheme and then use the equivariant
obstruction theory.

\subsection{The Target extension scheme}

\noindent According to the proposition \ref{prop:VezaProblem-Ekvivarijantan}
we would prefer to prove that there are no equivariant $\mathbb{Q}_{4n}$%
-maps $S^{3}\rightarrow T$. Thus the following scheme can be of some help.

\noindent \textbf{The basic idea} of the target extension scheme is to
\textit{find a }$\mathbb{Q}_{4n}$\textit{\ space }$E$\textit{\ which
contains the target space }$T$\textit{\ and to prove that there is no }$%
\mathbb{Q}_{4n}$\textit{-map }$S^{3}\rightarrow E$. Thus, this would imply
that there is no $\mathbb{Q}_{4n}$\textit{-map }$S^{3}\rightarrow T$.

\noindent Since the target space $T$ is the complement of the arrangement,
the basic idea can be refined as follows:

(A) \textit{Increase the codimension of the arrangement:} Take an arbitrary
hyper arrangement $\mathcal{J}$ in $W_{n}$. By the \textit{hyper}
arrangement we mean the arrangement of hyperplanes and / or closed
hyperplane halfspaces. Then form the minimal $\mathbb{Q}_{4n}$-invariant
arrangement $\mathcal{A}(\mathcal{J},\alpha )$ containing the family $%
\mathcal{J}\cap L(\alpha )=\{J\cap L(\alpha )~|~J\in \mathcal{J}\}$. Then
inclusion $\cup \mathcal{A}(\mathcal{J},\alpha )\subseteq \cup \mathcal{A}%
(\alpha )$ implies that $W_{n}-\cup \mathcal{A}(\mathcal{J},\alpha
)\supseteq W_{n}-\cup \mathcal{A}(\alpha )$. Observe that the dimension of
maximal elements of the arrangement $\mathcal{A}(\mathcal{J},\alpha )$ is $%
n-4$. Let us denote (when it suits us) the new union $\cup \mathcal{A}(%
\mathcal{J},\alpha )$ by $U^{\ast }$, and the new complement $W_{n}-\cup
\mathcal{A}(\mathcal{J},\alpha )$ by $T^{\ast }$.

(B) \textit{Apply the obstruction theory to the new question:} Is there a $%
\mathbb{Q}_{4n}$-map $S^{3}\rightarrow T^{\ast }$. Since the codimension of
maximal elements of the new defined arrangement in $W_{n}$ is $3$, the
target space $T^{\ast }$ is $1$-connected and consequently $2$-simple in the
sense that $\pi _{1}(T^{\ast })$ acts trivially on $\pi _{2}(T^{\ast })$.
Thus by Hurewicz theorem $\pi _{2}(T^{\ast })\cong \lbrack S^{2},T^{\ast
}]\cong H_{2}(T^{\ast };\mathbb{Z})$. Then the part of the obstruction exact
sequence (\cite{Dieck87}, \cite{guide2}) we are interested in is

\begin{equation*}
\begin{array}{ccccc}
\lbrack S^{3},T^{\ast }]_{\mathbb{Q}_{4n}} & \overset{\theta }{%
\longrightarrow } & \mathrm{Im}\left\{ [S_{(2)}^{3},T^{\ast }]_{\mathbb{Q}%
_{4n}}\longrightarrow \lbrack S_{(1)}^{3},T^{\ast }]_{\mathbb{Q}%
_{4n}}\right\} & \overset{\tau _{\mathbb{Q}_{4n}}}{\longrightarrow } & H_{%
\mathbb{Q}_{4n}}^{3}(S^{3},H_{2}(T^{\ast };\mathbb{Z}))%
\end{array}%
\end{equation*}%
where $S_{(1)}^{3}$ and $S_{(2)}^{3}$ are respectively the $1$- and $2$%
-skeleton of $S^{3}=P_{2n}\ast P_{2n}$ and $H_{2}(T^{\ast };\mathbb{Z})$ is
viewed as a $\mathbb{Q}_{4n}$-module. Since $[S_{(1)}^{3},T^{\ast }]_{%
\mathbb{Q}_{4n}}=\{\ast \}$ is a one-element set and $[S_{(2)}^{3},T^{\ast
}]_{\mathbb{Q}_{4n}}\neq \varnothing $, the sequence becomes
\begin{equation*}
\begin{array}{ccccc}
\lbrack S^{3},T^{\ast }]_{\mathbb{Q}_{4n}} & \longrightarrow & \{\ast \} &
\overset{\tau _{\mathbb{Q}_{4n}}}{\longrightarrow } & H_{\mathbb{Q}%
_{4n}}^{3}(S^{3},H_{2}(T^{\ast };\mathbb{Z})).%
\end{array}%
\end{equation*}%
The exactness means that the set $[S^{3},T^{\ast }]_{\mathbb{Q}_{4n}}\neq
\emptyset $ if and only if $\tau _{\mathbb{Q}_{4n}}(\ast )\in H_{\mathbb{Q}%
_{4n}}^{3}(S^{3},H_{2}(T^{\ast };\mathbb{Z}))$ is equal to zero. The element
$\tau _{\mathbb{Q}_{4n}}(\ast )$ depends only on $T^{\ast }$. Thus, the main
question transforms in

\subsection{A map in the general position}

\noindent We evaluate the class $\tau _{\mathbb{Q}_{4n}}(\ast )$ by the so
called \textquotedblright map in the general position\textquotedblright\
standard procedure, \cite{guide2}. Since the set $[S_{(2)}^{3},T^{\ast }]_{%
\mathbb{Q}_{4n}}$ is a point, it suffices to calculate an obstruction
cocycle for the particular map $h$.

\noindent Let $h:S^{3}\rightarrow W_{n}$ be an arbitrary $\mathbb{Q}_{4n}$%
-simplicial map which is in the general position. What we mean by the
general position is that for any simplex $\sigma $ in $S^{3}$
\begin{equation*}
h(\sigma )\cap U^{\ast }\neq \emptyset \text{ }\Rightarrow \mathrm{dim}%
(h(\sigma ))=3\text{ and }h(\sigma )\cap U^{\ast }=\{p_{1},..,p_{k}\}\subset
\mathrm{int}h(\sigma ).
\end{equation*}%
Now let $h:S^{3}\rightarrow W_{n}$ be a $\mathbb{Q}_{4n}$-map in the general
position. Then the associated cohomology class of the obstruction cocycle $%
c_{\mathbb{Q}_{4n}}(h)\in C_{\mathbb{Q}_{4n}}^{3}(S^{3},H_{2}(T^{\ast };%
\mathbb{Z}))=\mathrm{Hom}_{\mathbb{Q}_{4n}}(C_{3}(S^{3}),H_{2}(T^{\ast };%
\mathbb{Z}))$ for the map $h$ is equal to $\tau (\ast )$, i.e. $[c_{\mathbb{Q%
}_{4n}}(h)]=\tau (\ast )$. In order to calculate it let us describe the
cocycle $c_{\mathbb{Q}_{4n}}(h)$ a little bit closer. Let $\sigma $ be an
oriented $3$-simplex in $S^{3}$. Then $c_{\mathbb{Q}_{4n}}(h)(\sigma )\in
H_{2}(T^{\ast };\mathbb{Z})$ is the $h_{\ast }$ image of the fundamental
class of $\partial (\sigma )\cong S^{2}$ by the map $h_{\ast
}:H_{2}(\partial (\sigma );\mathbb{Z})\rightarrow H_{2}(T^{\ast };\mathbb{Z}%
) $, i.e.
\begin{equation*}
c_{\mathbb{Q}_{4n}}(h)(\sigma )=h_{\ast }[\partial (\sigma )].
\end{equation*}

\subsection{The nature of the obstruction cocycle}

\noindent The following proposition will narrow our attention to
the torsion part of the group \\
$H_{\mathbb{Q}_{4n}}^{3}(S^{3},H_{2}(T^{\ast };\mathbb{Z}))$ and
will be the perfect control factor in our calculations. It can be
also found in \cite{Bl-Vr-Ziv}.

\begin{proposition}
\label{prop:ObCoCycleTorzioni} The cohomology class of the obstruction
cocycle $c_{\mathbb{Q}_{4n}}(h)$ is a torsion element of the group $H_{%
\mathbb{Q}_{4n}}^{3}(S^{3},H_{2}(T^{\ast };\mathbb{Z}))$.
\end{proposition}

\begin{proof}
Let $H$ be a subgroup of $\mathbb{Q}_{4n}$. The restriction map

\begin{equation*}
r:H_{\mathbb{Q}_{4n}}^{3}(S^{3},H_{2}(T^{\ast };\mathbb{Z}))\rightarrow
H_{H}^{3}(S^{3},H_{2}(T^{\ast };\mathbb{Z}))
\end{equation*}%
on the cochain level sends $c\in C_{\mathbb{Q}_{4n}}^{3}(S^{3},H_{2}(T^{\ast
};\mathbb{Z}))$ a $\mathbb{Q}_{4n}$-cochain to now $H$-cochain $c\in
C_{H}^{3}(S^{3},H_{2}(T^{\ast };\mathbb{Z}))$. The definition of the
obstruction cocycle implies that $r(c_{\mathbb{Q}_{4n}}(h))$ is the
obstruction cocycle for the extension of the $H$-map in the general position
$h$. It is a known fact (\cite{Brown} Section III.9. Proposition 9.5.(ii))
that the composition of the restriction with the transfer $\tau
:H_{H}^{3}(S^{3},H_{2}(T^{\ast };\mathbb{Z}))\rightarrow H_{\mathbb{Q}%
_{4n}}^{3}(S^{3},H_{2}(T^{\ast };\mathbb{Z}))$ is just a multiplication by
the index $[\mathbb{Q}_{4n}:H]$,
\begin{equation*}
\begin{array}{ccccc}
H_{\mathbb{Q}_{4n}}^{3}(S^{3},H_{2}(T^{\ast };\mathbb{Z})) & \rightarrow &
H_{H}^{3}(S^{3},H_{2}(T^{\ast };\mathbb{Z})) & \rightarrow & H_{\mathbb{Q}%
_{4n}}^{3}(S^{3},H_{2}(T^{\ast };\mathbb{Z})) \\
\lbrack c_{\mathbb{Q}_{4n}}(h)] & \longmapsto & [c_{H}(h)] & \longmapsto & [%
\mathbb{Q}_{4n}:H]\cdot \lbrack c_{\mathbb{Q}_{4n}}(h)]\text{.}%
\end{array}%
\end{equation*}%
Particularly, let $H$ be the trivial subgroup of $\mathbb{Q}_{4n}$. The
sphere $S^{3}$ is $2$-connected and $H$ is\ a trivial group, so the map $h$
can be extended to a $H$-map $S^{3}\rightarrow M$. Thus, $[c_{H}(h)]=0$ and
consequently $[\mathbb{Q}_{4n}:H]\cdot \lbrack c_{\mathbb{Q}_{4n}}(h)]=0$ in
$H_{\mathbb{Q}_{4n}}^{3}(S^{3},H_{2}(T^{\ast };\mathbb{Z}))$.
\end{proof}

\subsection{The $\mathbb{Q}_{4n}$ cellular structures on $S^{3}$}

\noindent In order to start efficient computations of the obstruction
cocycle we need to describe the concrete $\mathbb{Q}_{4n}$ $CW$-structures
of the sphere $S^{3}$. The proposition \ref{prop:promenaDejstva} allows us
to be very picky in selecting the adequate cellular structures. We describe
two $\mathbb{Q}_{4n}$ cellular structures the natural one and the most
economical one, and the cellular map joining them. The direct consequence of
these discussions is the isomorphism
\begin{equation*}
H_{\mathbb{Q}_{4n}}^{3}(S^{3},H_{2}(T^{\ast };\mathbb{Z}))\cong
H_{2}(T^{\ast };\mathbb{Z})_{\mathbb{Q}_{4n}}.
\end{equation*}

\begin{figure}[htb]
\centering
\includegraphics[scale=0.70]{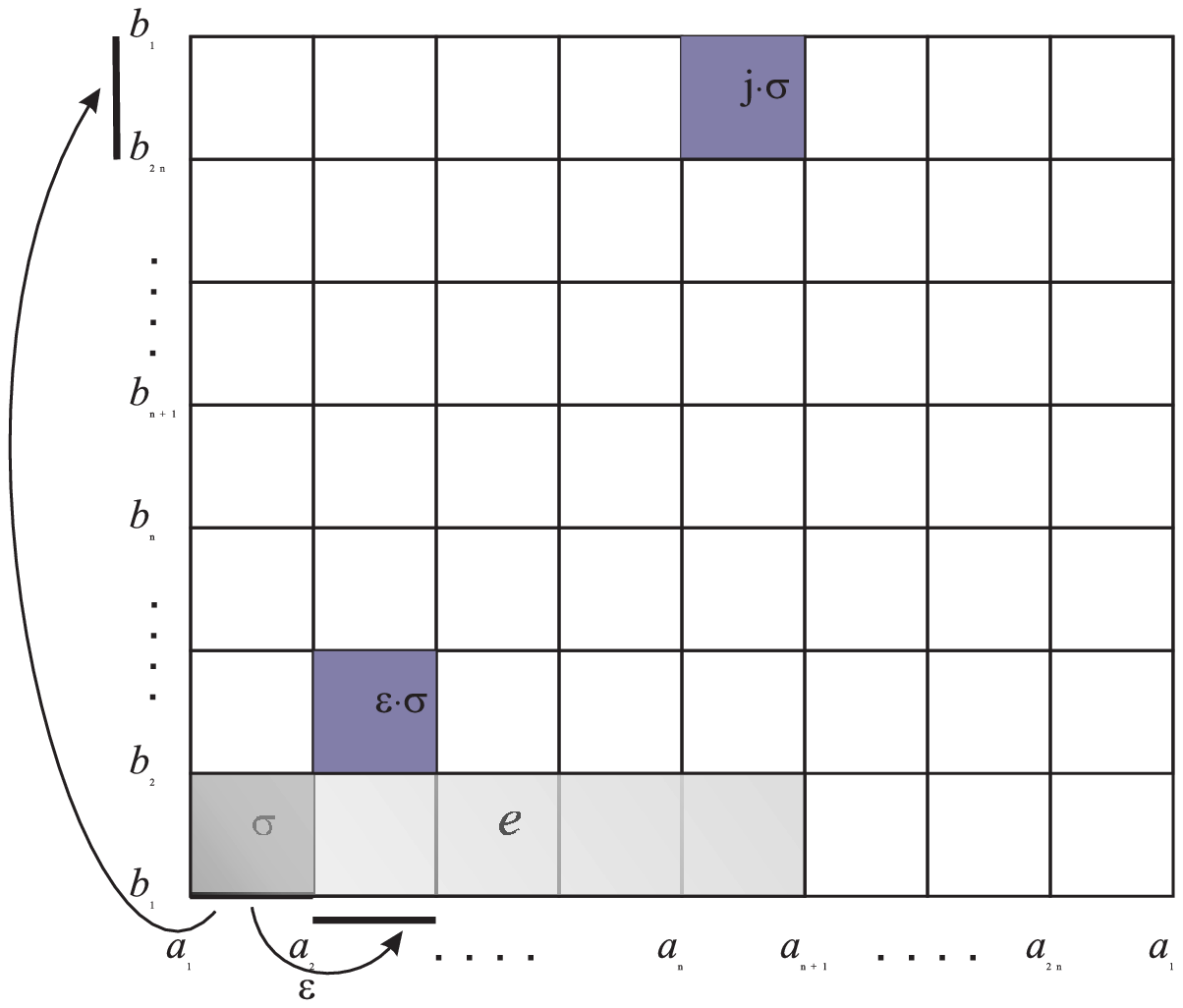}
\caption{The model for two $\mathbf{Q}_{4n}$ cellular structures
of $S^3$.} \label{fig:Fig1}
\end{figure}

\noindent \textbf{The natural }$\mathbb{Q}_{4n}$ \textbf{cellular-simplicial
structure.} This structure comes from the join decomposition $%
S^{3}=S^{1}\ast S^{1}$ of the $3$-sphere. Let the sphere $S^{1}$ be
represented by the simplicial complex of the regular $2n$-gon $P_{2n}$. Then
the sphere $S^{3}$, as the simplicial complex, is the join $P_{2n}^{(1)}\ast
P_{2n}^{(2)}$ of two copies of $P_{2n}$. Let the vertex of $P_{2n}^{(1)}$
and $P_{2n}^{(2)}$ be denoted by $a_{1},..,a_{2n}$ and $b_{1},..,b_{2n}$,
respectively. The action of the group $\mathbb{Q}_{4n}$ on $S^{3}$ is
defined on vertices by%
\begin{equation*}
\epsilon \cdot a_{i}=a_{i\ \textrm{mod}\ 2n+1}\text{, }\epsilon \cdot b_{i}=b_{i%
\ \textrm{mod}\ 2n+1}\text{, }j\cdot a_{1}=b_{1}
\end{equation*}%
and it extends equivariantly to upper skeletons. Then, for example%
\begin{equation*}
\begin{array}{c}
\text{ }j\cdot a_{i}=j\epsilon ^{i-1}\cdot a_{1}=\epsilon
^{(2n-1)(i-1)}j\cdot a_{1}=\epsilon ^{(2n-1)(i-1)}\cdot b_{1}=b_{(2n-i+1)%
\ \textrm{mod}\ 2n+1} \\
j\cdot \lbrack a_{1},a_{2};b_{1},b_{2}]=[b_{1},b_{2n};a_{n+1},a_{n}].%
\end{array}%
\end{equation*}%
The associated chain complex $\mathfrak{C}=\{\emph{C}_{i}\}$ has the form%
\begin{equation*}
0\longrightarrow \mathbb{Z}^{4n^{2}}\longrightarrow \mathbb{Z}%
^{8n^{2}}\longrightarrow \mathbb{Z}^{4n^{2}+4n}\longrightarrow \mathbb{Z}%
^{4n}\longrightarrow 0.
\end{equation*}%
This $\mathbb{Q}_{4n}$ cellular structure makes the beautiful
rectangular middle section of the join representation
$S^{3}=S^{1}\ast S^{1}=[0,1]\ast \lbrack 0,1]/\approx $, (with
some identifications). The action is indicated in the figure
\ref{fig:Fig1}.


\noindent \textbf{The economic }$\mathbb{Q}_{4n}$ \textbf{cellular structure.%
} The second one comes from the minimal resolution of $\mathbb{Z}$ by free $%
\mathbb{Q}_{4n}$-modules described in \cite{CaEi}, p.~253. The associated
cellular complex has one $\mathbb{Q}_{4n}$ $0$-cell $a$, two $\mathbb{Q}_{4n}
$ $1$-cells $b$ and $b^{\prime }$, two $\mathbb{Q}_{4n}$ $2$-cells $c$ and $%
c^{\prime }$, and finally then one $\mathbb{Q}_{4n}$ $3$-cell $e$. The
associated chain complex $\mathfrak{D}=\{\emph{D}_{i}\}$ has the form
\begin{equation*}
0\rightarrow \mathbb{Z}[\mathbb{Q}_{4n}]e\overset{\partial }{\rightarrow }%
\mathbb{Z}[\mathbb{Q}_{4n}]c\oplus \mathbb{Z}[\mathbb{Q}_{4n}]c^{\prime }%
\overset{\partial }{\rightarrow }\mathbb{Z}[\mathbb{Q}_{4n}]b\oplus \mathbb{Z%
}[\mathbb{Q}_{4n}]b^{\prime }\overset{\partial }{\rightarrow }\mathbb{Z}[%
\mathbb{Q}_{4n}]e\rightarrow 0
\end{equation*}%
where
\begin{equation*}
\begin{array}{lll}
\partial e=(\epsilon -1)c-(\epsilon j-1)c^{\prime } & \partial
c=(1+..+\epsilon ^{n-1})b-(j+1)b^{\prime } & \partial c^{\prime }=(\epsilon
j+1)b+(\epsilon -1)b^{\prime } \\
& \partial b=(\epsilon -1)a & \partial b^{\prime }=(j-1)a%
\end{array}%
.
\end{equation*}%
Thus, it will be enough to look at the obstruction cocycle $c_{\mathbb{Q}%
_{4n}}(h)$ on the maximal cell $e$, and to prove that its image is or is not
zero, when we pass to cohomology.

\noindent \textbf{Computing equivariant cohomology.} The explicit formulas
for the chain complex $\{\emph{D}_{i}\}$ is lurking us to apply the $\mathrm{%
Hom}_{\mathbb{Q}_{4n}}(\cdot ,H_{2}(T^{\ast };\mathbb{Z}))$ functor. The
result is the equivariant cochain complex
\begin{equation*}
0\longleftarrow H_{2}(T^{\ast };\mathbb{Z})\overset{\Gamma }{\longleftarrow }%
H_{2}(T^{\ast };\mathbb{Z})\oplus H_{2}(T^{\ast };\mathbb{Z})\leftarrow
H_{2}(T^{\ast };\mathbb{Z})\oplus H_{2}(T^{\ast };\mathbb{Z})\longleftarrow
H_{2}(T^{\ast };\mathbb{Z})\longleftarrow 0
\end{equation*}%
where $\Gamma (p,q)=(\epsilon -1)p-(\epsilon j-1)q$ for $p,q\in
H_{2}(T^{\ast };\mathbb{Z})$. The definition of the equivariant cohomology
and a standard calculation imply that
\begin{equation*}
H_{\mathbb{Q}_{4n}}^{3}(S^{3},H_{2}(T^{\ast };\mathbb{Z}))=H_{2}(T^{\ast };%
\mathbb{Z})/\mathrm{Im}\Gamma \cong H_{2}(T^{\ast };\mathbb{Z})_{\mathbb{Q}%
_{4n}}
\end{equation*}%
where $H_{2}(T^{\ast };\mathbb{Z})_{\mathbb{Q}_{4n}}$ denote the group of
coinvariants of the $\mathbb{Q}_{4n}$-module $H_{2}(T^{\ast };\mathbb{Z})$,
\cite{Brown}.

\noindent \textbf{The chain map.} There exist a cellular map $\mathfrak{f}:%
\mathfrak{D}\rightarrow \mathfrak{C}$, which can be visualized on the figure %
\ref{fig:Fig1} by identifying top dimensional cell $e$ of the
economic cell structure with the transparent-shaded fundamental
domain.

\noindent Since we are interested in $3$-cochains, i.e. elements of $C_{%
\mathbb{Q}_{4n}}^{3}(S^{3},H_{2}(T^{\ast };\mathbb{Z}))$, we only need to
know the concrete expression for the cellular map on the top dimensional
cell $e$,%
\begin{equation*}
\mathfrak{f}(e)=[a_{1},\epsilon a_{1}]\ast \lbrack b_{1},\epsilon
b_{1}]+[\epsilon a_{1},\epsilon ^{2}a_{1}]\ast \lbrack b_{1},\epsilon
b_{1}]+...+[\epsilon ^{n-1}a_{1},\epsilon ^{n}a_{1}]\ast \lbrack
b_{1},\epsilon b_{1}]\text{,}
\end{equation*}%
where the simplexes on the right hand side are appropriately oriented.

\subsection{\label{sec:pointClasses}Point classes}

\noindent It is of utmost importance for computation of the obstruction
cocycle $c(h)_{\mathbb{Q}_{4n}}$ to pinpoint some elements from the
coefficient module $H_{2}(T^{\ast };\mathbb{Z})$. The similar discussion can
be found in \cite{Bl-Vr-Ziv}.

\noindent Let $\{W_{1},W_{2},\ldots ,W_{k}\}$ be the family on maximal
elements of the arrangement $\mathcal{A}$ of linear (closed half-)subspaces
in an $(n+m)$-dimensional, Euclidean space $E$. Let also maximal elements
have the constant dimension $\dim W_{i}=n$. Let $\hat{D}(\mathcal{A})=\cup
\mathcal{A}\cup \{+\infty \}\subset E\cup \{+\infty \}\cong S^{n+m}$ be the
compactified union and $M(\mathcal{A})=E\setminus \cup \mathcal{A}$ the
complement of the arrangement.

\noindent Let the point $x\in \mathrm{int}(W_{i}\setminus \cup _{j\neq
i}W_{j})$ and let $D_{\varepsilon }(x)=x+D_{\varepsilon }$ be a disc around $%
x$, where $D_{\varepsilon }$ is the $\varepsilon $-disc in the orthogonal
complement $W_{i}^{\perp }$. For a sufficiently small $\varepsilon $ the
intersection $D_{\epsilon }(x)\cap (\cup _{j\neq i}W_{j})$ vanishes and we
assume fix such an $\varepsilon $. We assume that $D_{\epsilon }(x)$ is
oriented by the orientation inherited from the ambient orientation and the
orientation prescribed of $W_{i}$.

\noindent \textbf{The point class} $[x]\in H_{m}(E,M(\mathcal{A});\mathbb{Z}%
) $ of $x$ is inclusion image of the fundamental class of the pair $%
(D_{\epsilon }(x),\partial D_{\epsilon }(x))$. By the Excision axiom
described $\varepsilon $ has no effect on the class $[x]$. In addition, by
the Homotopy axiom, $[x]$ is uniquely determined by the connected component
of $W_{i}\setminus \cup _{j\neq i}W_{j}$. The image $\left\Vert x\right\Vert
:=\partial \lbrack x]$ of the point class $[x]$ by the isomorphism $%
H_{m}(E,M(\mathcal{A});\mathbb{Z})\rightarrow H_{m-1}(M(\mathcal{A});\mathbb{%
Z})$ is also called the \textit{point class }of\emph{\ }$x$ and has all the
properties of the original one.


\begin{figure}[htb]
\centering
\includegraphics[scale=1]{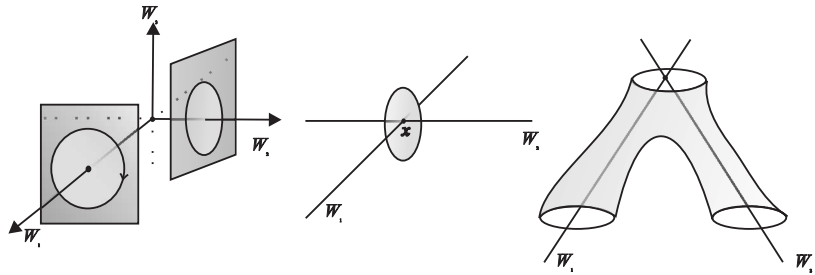}
\caption{The point classes and the broken point classes}
\label{fig:Fig2}
\end{figure}

\begin{proposition}
\label{prop:PointClass}Let $\mathcal{A}$ be the arrangement of linear
subspaces as above. Let $x\in \mathrm{int}(W_{i}\setminus \cup _{j\neq
i}W_{j})=W_{i}\setminus \cup _{j\neq i}W_{j}$.

(A) The class $[x]\in H_{m}(E,M(\mathcal{A});\mathbb{Z})$ does not vanish.

(B) If we assume that $(\forall j)~i\neq j~\Rightarrow ~\mathrm{codim}%
_{W_{i}}(W_{i}\cap W_{j})>1$, then the class $[x]$ does not depend on
particular $x$.

(C) If there is a subspace $W_{j}$ such that $\mathrm{codim}%
_{W_{i}}(W_{i}\cap W_{j})=1$ and $x_{1}$, $x_{2}$ belong to different
connected components of $W_{i}\setminus \cup _{j\neq i}W_{j}$, then $%
[x_{1}]\neq \lbrack x_{2}]$.
\end{proposition}

\begin{proof}
(A) By the Ziegler-\v{Z}ivaljevi\'{c} formula \cite{ZZ}, has he wedge
decomposition%
\begin{equation*}
\hat{D}(\mathcal{A})\simeq \hat{W}_{1}\vee \hat{W}_{2}\vee \ldots \vee \hat{W%
}_{k}\vee ...
\end{equation*}%
where the displayed factors correspond to maximal elements. Let $\Delta (%
\hat{W}_{i})\in H^{m}(E,M(\mathcal{A});\mathbb{Z})$ be the Poincar\'{e}%
-Alexander dual of the fundamental homology class $[\hat{W}_{i}]\in H_{n}(%
\hat{D}(\mathcal{A});\mathbb{Z})$ associated with the sphere $\hat{W}_{i}$
in $\hat{E}$. Then $\Delta (\hat{W}_{i})([N])\in \mathbb{Z}$ is the
intersection number $[\hat{W}_{i}]\cap \lbrack N]$, whenever this number is
correctly defined, for example if $N$ is a manifold and the intersection $%
\hat{W}_{i}\cap N$ is transversal. From here we see that $\Delta (\hat{W}%
_{i})([x])=1$ and so $[x]$ does not vanish.

(B) Since the assumption implies the connectness of the complement $%
W_{i}\setminus \cup _{j\neq i}W_{j}$, the statement follows by some homotopy.

(C) When $\mathrm{codim}_{W_{i}}(W_{i}\cap W_{j})=1$ the subspaces $W_{i}$
and $W_{j}$ are decomposed by the hyperplane $W_{i}\cap W_{j}$ into unions
of closed half-spaces, $W_{i}=W_{i}^{1}\cup W_{i}^{2}$ and $%
W_{j}=W_{j}^{1}\cup W_{j}^{2}$, respectively. Let us assume that $x_{1}\in
W_{i}^{1}$ and $x_{2}\in W_{i}^{2}$. The half-space $W_{i}^{1}$ and $%
W_{j}^{1}$ glued along the common boundary and compactified determine a
sphere $S\subset \hat{E}$. Since the Ziegler-\v{Z}ivaljevi\'{c}
decomposition \cite{ZZ} involves a choice of generic points, the points can
be chosen in such a way that the sphere $S$ appears in the decomposition.
Precisely, the sphere $S$ is the factor $\widehat{W_{i}\cap W_{j}}\ast
\Delta (P_{<W_{i}\cap W_{j}})\cong S^{n-1}\ast S^{0}\cong S^{n}$, where $P$
denote the intersection poset of the arrangement $\mathcal{A}$. thus, the
fundamental class $[S]$ is nontrivial in $H_{n}(\hat{D}(\mathcal{A});\mathbb{%
Z})$. Like in the previous case, let $\Delta (S)\in H^{m}(E,M(\mathcal{A});%
\mathbb{Z})$ be the Poincar\'{e}-Alexander dual to $[S]$. Then%
\begin{equation*}
\Delta (S)([x_{1}])=\pm 1\text{ and }\Delta (S)([x_{2}])=0
\end{equation*}%
and consequently $[x_{1}]\neq \lbrack x_{2}]$.
\end{proof}

\noindent \textbf{The broken point class.} In order to simplify the
exposition let $\mathcal{A}$ be a very special arrangement consisting of two
vector spaces $W_{1}$ and $W_{2}$ of the same dimension $m$ inside the
vector space $E=\mathbb{R}^{n+m}$. Let $U=W_{1}\cap W_{2}$ be the
intersection and $k=\dim (W_{1}\cap W_{2})$. Let us also fix $V$, a vector
space of dimension $n$ such that $E=W_{1}\oplus V=W_{2}\oplus V$, where $%
\oplus $ denotes the direct sum of vector spaces. For $x\in W_{1}\cap W_{2}$
let $D_{\varepsilon }(x)$ be a disc $x+D_{\varepsilon }$, where $%
D_{\varepsilon }$ is the $\varepsilon $-disc in $V$. Then the \textit{broken
point class} of $x$, denoted by $[x]\in H_{n}(E,E\backslash W_{1}\cup W_{2};%
\mathbb{Z})$, is the fundamental class of the pair $(D_{\epsilon
}(x),\partial D_{\epsilon }(x))$. In addition the image $\left\Vert
x\right\Vert :=\partial \lbrack x]$ of the broken point class $[x]$ by the
isomorphism $H_{n}(E,E\backslash W_{1}\cup W_{2};\mathbb{Z})\rightarrow
H_{n-1}(E\backslash W_{1}\cup W_{2};\mathbb{Z})$ is also called the broken
point class\textit{\ }of\emph{\ }$x$ and has all the properties of the
original one.

\noindent An illustration for the case $m=1$ and $n=2$ can be seen
in the figure \ref{fig:Fig2}. The same definition stands even if
$W_{1}$ and $W_{2}$ are closed half spaces intersecting over
linear space $U$. The broken point classes can be naturally
expressed as linear combinations of the ordinary point classes.
The third picture in the figure \ref{fig:Fig2} illustrates this
situation.


\begin{figure}[htb]
\centering
\includegraphics[scale=1.3]{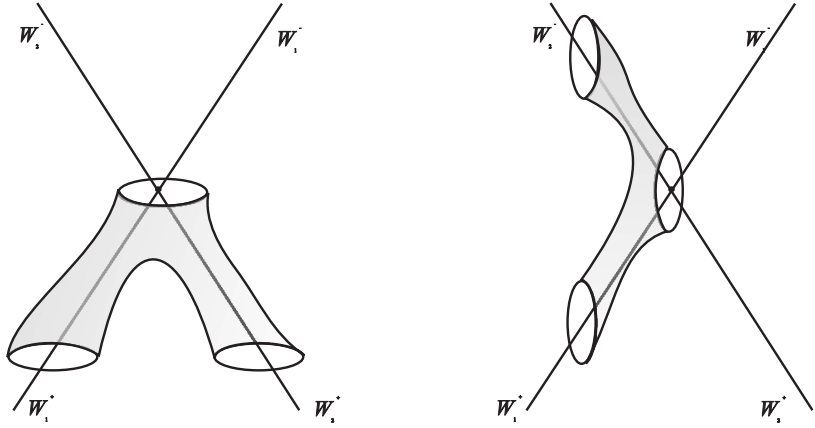}
\caption{Two posible point classes in the case $dim(W1 \cap
W2)=m-1$} \label{fig:Fig5}
\end{figure}

\noindent Let us now suppose that $k=\dim (W_{1}\cap W_{2})=\dim
(W_{1})-1=m-1$. Then there are two different broken point classes
associated to every point $x\in W_{1}\cap W_{2}$, i.e. the broken
point class $[x]$ depends on the choice of the vector space $V$.
Let us omit the format proof of this fact and illustrate this
situation by the figure \ref{fig:Fig5}. The
formal proof would go along the lines of the proposition \ref%
{prop:PointClass}. The main property that differs these two classes,
brusquely explained, is that if we substitute each sphere with the disc and
move these discs a little bit:

(A) the first disc will intersect simultaneously $W_{1}^{+}$ and $W_{2}^{+}$%
, or $W_{1}^{-}$ and $W_{2}^{-}$; and

(B) the second disc will intersect simultaneously $W_{1}^{+}$ and $W_{2}^{-}$%
, or $W_{1}^{-}$ and $W_{2}^{+}$.

\subsection{Calculating $H_{2}(T^{\ast };\mathbb{Z})$ and $H_{2}(T^{\ast };%
\mathbb{Z})_{\mathbb{Q}_{4n}}$}

\noindent \textbf{The homology. }In order to analyze the second homology of
the complement $T^{\ast }$ we use the following equality%
\begin{equation*}
W_{n}\backslash \cup \mathcal{A}(\mathcal{J},\alpha )=S^{n-1}\backslash \cup
\widehat{\mathcal{A}}(\mathcal{J},\alpha )
\end{equation*}%
where $\widehat{\mathcal{A}}(\mathcal{J},\alpha )$ denotes the one point
compactification of the arrangement $\mathcal{A}(\mathcal{J},\alpha )$. This
allow us to use the Poincar\'{e}-Alexander duality and to work with the
arrangement $\widehat{\mathcal{A}}(\mathcal{J},\alpha )$ instead of its
complement. The Poincar\'{e}-Alexander duality and the Universal Coefficient
isomorphisms give us the sequence of isomorphisms (assuming $\mathbb{Z}$
coefficients)
\begin{eqnarray}
H_{2}(W_{n}\backslash \cup \mathcal{A}(\mathcal{J},\alpha ))
&=&H_{2}(S^{n-1}\backslash \cup \widehat{\mathcal{A}}(\mathcal{J},\alpha
))\cong H^{(n-1)-2-1}(\cup \widehat{\mathcal{A}}(\mathcal{J},\alpha ))
\label{Poincare-Alexander 1} \\
&\cong &\mathrm{Hom}(H_{n-4}(\cup \widehat{\mathcal{A}}(\mathcal{J},\alpha
)),\mathbb{Z})\oplus \mathrm{Ext}(H_{n-3}(\cup \widehat{\mathcal{A}}(%
\mathcal{J},\alpha )),\mathbb{Z}).  \notag
\end{eqnarray}%
Since the maximal elements of the arrangement $\mathcal{A}(\mathcal{J}%
,\alpha )$ are $(n-4)$-dimensional linear (half-) subspaces, respectively,
the $\mathrm{Ext}$ factor in the above equality vanish. Precisely, $%
H_{n-3}(\cup \widehat{\mathcal{A}}(\mathcal{J},\alpha );\mathbb{Z})=0$.
Therefore,%
\begin{equation}
H_{2}(W_{n}\backslash \cup \mathcal{A}(\mathcal{J},\alpha );\mathbb{Z})\cong
\mathrm{Hom}(H_{n-4}(\cup \widehat{\mathcal{A}}(\mathcal{J},\alpha );\mathbb{%
Z}),\mathbb{Z}).  \label{Isomorphisms Complement - Arrangement}
\end{equation}%
The Ziegler-\v{Z}ivaljevi\'{c} formula implies the following homology
decomposition (assuming $\mathbb{Z}$ coefficients)
\begin{equation}
H_{n-4}(\cup \widehat{\mathcal{A}}(\mathcal{J},\alpha );\mathbb{Z})\cong
\underset{d=0}{\overset{n-4}{\bigoplus }}\underset{V\in P(\alpha ):\dim V=d}%
{\bigoplus }H_{n-4}(\Delta (P_{<V})\ast \hat{V}) \label{Z-Z
homology decomposition - 1}
\end{equation}%
where $P(\alpha )$ is the intersection poset of the arrangement $\mathcal{A}%
(J,\alpha )$ and $\hat{V}$ one-point compactification of the element $V$.
Thus, the property of the hom functor imply%
\begin{equation}
H_{2}(W_{n}\backslash \cup \mathcal{A}(\mathcal{J},\alpha
);\mathbb{Z})\cong \underset{d=0}{\overset{n-4}{\bigoplus
}}\mathrm{Hom}\left( \underset{V\in
P(\alpha ):\dim V=d}{\bigoplus }H_{n-4}(\Delta (P_{<V})\ast \hat{V});%
\mathbb{Z}\right)  \label{hom decomposition - 1}
\end{equation}%
\noindent \textbf{Poincar\'{e} dual of the (broken) point class.} Let us
illustrate the isomorphism (\ref{Isomorphisms Complement - Arrangement}) on
the point classes for the general arrangement of linear (half-) spaces $%
\mathcal{A}$ in $\mathbb{R}^{n+m}$. Let maximal elements $%
\{W_{1},W_{2},\ldots ,W_{k}\}$ have the constant dimension $\dim W_{i}=n$.
Let $x\in W_{1}\backslash \cup _{i\neq 1}W_{i}$ and $S=\partial
D_{\varepsilon }(x)$, where $D_{\varepsilon }(x)=x+D_{\varepsilon }$ and $%
D_{\varepsilon }$ is a small disk in $W_{1}{}^{\bot }$. Then the isomorphism
, $\vartheta :H_{m-1}(\mathbb{R}^{n+m}\backslash \cup \mathcal{A};\mathbb{Z}%
)\rightarrow \mathrm{Hom}(H_{n}(\cup \widehat{\mathcal{A}};\mathbb{Z}),%
\mathbb{Z})$ can be expressed for $t\in H_{n}(\cup \widehat{\mathcal{A}};%
\mathbb{Z});\mathbb{Z})$ by%
\begin{equation}
\vartheta (\left\Vert x\right\Vert ):H_{n}(\cup \widehat{\mathcal{A}};%
\mathbb{Z})\rightarrow \mathbb{Z}\text{, }\vartheta (\left\Vert x\right\Vert
)(t)=\mathrm{link}(S,T)  \label{link-isomorphism}
\end{equation}%
if $T$ is a submanifold in $\cup \widehat{\mathcal{A}}$ representing the
homology class $t$ and linking number $\mathrm{link}(S,T)$ is correctly
defined.

\begin{example}
Let the arrangement $\mathcal{A}$ be given by the figure
\ref{fig:Fig4},(A). Then the Hasse diagram of the arrangement
$\mathcal{A}$ is like the in the figure \ref{fig:Fig4},(B).

\begin{figure}[htb]
\centering
\includegraphics[scale=1]{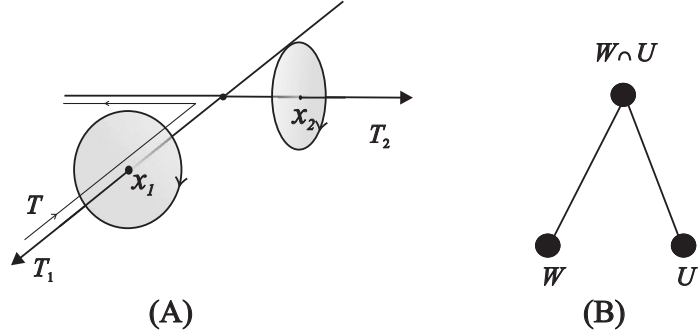}
\caption{The arrangement $\mathcal{A}$ with indicated generators
in homology.} \label{fig:Fig4}
\end{figure}

Then $H_{1}(\cup \widehat{\mathcal{A}};\mathbb{Z})\cong \mathbb{%
Z\oplus Z\oplus Z}$, where the generators / spheres $T$, $T_{1}$ and $T_{2}$
indicated in the figure. Then we can read of isomorphism $\vartheta :H_{1}(%
\mathbb{R}^{3}\backslash \cup \mathcal{A};\mathbb{Z})\rightarrow \mathrm{Hom}%
(H_{1}(\cup \widehat{\mathcal{A}};\mathbb{Z}),\mathbb{Z})$ for point classes
$\left\Vert x_{1}\right\Vert $ and $\left\Vert x_{2}\right\Vert $ from the
picture,%
\begin{eqnarray*}
\vartheta (\left\Vert x_{1}\right\Vert )(t_{1}) &=&1,\vartheta (\left\Vert
x_{1}\right\Vert )(t_{2})=0,\vartheta (\left\Vert x_{1}\right\Vert )(t)=-1%
\text{,} \\
\vartheta (\left\Vert x_{2}\right\Vert )(t_{1}) &=&0,\vartheta (\left\Vert
x_{2}\right\Vert )(t_{2})=1,\vartheta (\left\Vert x_{2}\right\Vert )(t)=0%
\text{.}
\end{eqnarray*}
\end{example}

\noindent \textbf{The coinvariants. }When time comes to compute coinvariants
the best news would be that the Poincar\'{e}-Alexander duality map and the
Universal coefficient isomorphism are equivariant maps. Then the
isomorphisms (\ref{Isomorphisms Complement - Arrangement}) would be
isomorphisms of $\mathbb{Q}_{4n}$-modules, and there would be no difference
in what module we work. (Un)fortunately, the Poincar\'{e}-Alexander duality
map is not a $\mathbb{Q}_{4n}$-map, but a $\mathbb{Q}_{4n}$-map up to a
orientation character, while Universal coefficient isomorphism is an $%
\mathbb{Q}_{4n}$-map. If $o$ is an orientation of the $\mathbb{Q}_{4n}$%
-sphere $S^{n-1}$ or $S^{2(n-1)}$, then $o$ determines Poincar\'{e}%
-Alexander duality map $\gamma _{o}$ \cite{Mu}. In particular, $g\cdot
o=\det (g)\cdot o$, where $g\in \mathbb{Q}_{4n}\subseteq GL_{n}(\mathbb{R})$%
, or $g\in \mathbb{Q}_{4n}\subseteq GL_{2n}(\mathbb{R})$.

\noindent The natural $\mathbb{Q}_{4n}$-action on the union $\cup \widehat{%
\mathcal{A}}(\mathcal{J},\alpha )$, inherited from the ambient $\mathbb{Q}%
_{4n}$-action, respects the dimensional decomposition (\ref{Z-Z homology
decomposition - 1}) of the $(n-4)$-homology. When these homologies are free,
then the decomposition (\ref{hom decomposition - 1}) becomes true even
without $\mathrm{Hom}$. Let us assume that $(n-4)$-homology of the
arrangement $\widehat{\mathcal{A}}(\mathcal{J},\alpha )$ is free. In some
calculation we will see that we actually do not need the whole group to be
free, but just some factors in the decomposition (\ref{hom decomposition - 1}%
). Therefore, if we would rather work with $\mathbb{Q}_{4n}$-module $%
H_{n-4}(\cup \widehat{\mathcal{A}}(\mathcal{J},\alpha );\mathbb{Z})$ instead
of $\mathbb{Q}_{4n}$-module $H_{2}(W_{n}\backslash \cup \mathcal{A}(\mathcal{%
J},\alpha );\mathbb{Z})$ we have to modify the $\mathbb{Q}_{4n}$-action.
Specifically, let $l\in H_{n-4}(\cup \widehat{\mathcal{A}}(\mathcal{J}%
,\alpha );\mathbb{Z})$ and $g\in \mathbb{Q}_{4n}$, then
\begin{equation*}
g\ast l=\det (g)~g\cdot l
\end{equation*}%
where $\ast $ is the new modified action, and $\cdot $ the old one. Let $%
\sim $ denote the congruence relation on $H_{n-4}(\cup \widehat{\mathcal{A}}(%
\mathcal{J},\alpha ),\mathbb{Z})$ which class of zero is a subgroup
generated by the elements of the form $g\ast x-x$, $g\in \mathbb{Q}_{4n}$, $%
x\in H_{n-4}(\widehat{\mathcal{A}}(\mathcal{J},\alpha ),\mathbb{Z})$. Then
there is an isomorphism%
\begin{equation}
H_{2}(W_{n}\backslash \cup \mathcal{A}(\mathcal{J},\alpha );\mathbb{Z})_{%
\mathbb{Q}_{4n}}\cong H_{n-4}(\cup \widehat{\mathcal{A}}(\mathcal{J},\alpha
),\mathbb{Z})/\sim   \label{isomorphism of the modified action}
\end{equation}

\noindent \textbf{The torsion of the coinvariants.} The proposition \ref%
{prop:ObCoCycleTorzioni} directs us to search for the obstruction cohomology
class in the torsion part of the coinvariants. Therefore, if the
coinvariants are free, then the obstruction is zero and the map exists. Thus
everythig we have done is in vain. It is of utmost importance to get a
feeling, when some $\mathbb{G}$ arrangement can produce nontrivial torsion
group in appropriate coinvariant group.

\noindent Let us discuss a few simple examples, which met all the above
assumptions. These examples come from the fan partition problems, and we
work with cyclic groups in order to simplify computations.

\begin{example}
\label{Ex:Coinvariants1}Let $\mathcal{A}$ be the minimal $\mathbb{Z}_{8}$
arrangement in $\mathbb{R}^{8}$ containing subspace
\begin{equation*}
L=\{\mathbf{x}\in \mathbb{R}%
^{4}~|~x_{1}+x_{2}=x_{3}+x_{4}=x_{5}+x_{6}=x_{7}+x_{8}=0\}\subset W_{8}.
\end{equation*}%
The group $\mathbb{Z}_{8}=\langle \varepsilon \rangle $ acts by a cyclic
permutation, i.e. $\varepsilon \cdot (x_{1},..,x_{8})=(x_{8},x_{1},..,x_{7})$%
. Then $\det (\varepsilon )=(-1)^{8+1}$.

\noindent It is not hard to see that $\mathcal{A}=\{L,\varepsilon L,L\cap
\varepsilon L\}$, $L\cap \varepsilon L=\{0\}$, and consequently the by Z-Z
decomposition (\ref{Z-Z homology decomposition - 1}),
\begin{equation*}
H_{4}(\cup \widehat{\mathcal{A}};\mathbb{Z})\cong \underset{d=0}{\overset{4}{%
\bigoplus }}\underset{V\in P:\dim V=d}{\bigoplus }H_{4}(\Delta
(P_{<V})\ast \hat{V};\mathbb{Z})\cong 0\oplus 0\oplus 0\oplus (\mathbb{Z}%
\oplus \mathbb{Z})
\end{equation*}%
where $P$ is the intersection poset of the arrangement $\mathcal{A}$. Since,
homology $H_{4}(\cup \widehat{\mathcal{A}};\mathbb{Z})$ has no torsion we
have that $H_{2}(W_{8}\backslash \cup \mathcal{A};\mathbb{Z})\cong \mathbb{Z}%
\oplus \mathbb{Z}$.

\noindent Let $l$ and $\varepsilon \cdot l$ be the generators of $H_{4}(\cup
\widehat{\mathcal{A}};\mathbb{Z})$ corresponding to spheres $\widehat{L}$
and $\widehat{\varepsilon L}$. The the equality $L=\varepsilon ^{2}L$ in $%
\mathbb{R}^{8}$ imply the equality $l=\epsilon (\varepsilon ^{2}\cdot l)$ in
$H_{4}(\cup \widehat{\mathcal{A}};\mathbb{Z})$, where $\epsilon \in \{1,-1\}$%
. The sign $\epsilon $ depends on the nature of $\varepsilon ^{2}$. If $%
\varepsilon ^{2}$ changes the orientation of $L$, then $\epsilon =-1$,
otherwise $\epsilon =1$. To calculate the sign $\epsilon $ we use a $%
\varepsilon ^{2}$-invariant decomposition $\mathbb{R}^{8}=L\oplus L^{\bot }$%
, where $L^{\bot }$ denote the orthogonal complement of $L$. Then,
\begin{equation*}
\det_{\mathbb{R}^{8}}\varepsilon ^{2}=\det_{L}\varepsilon ^{2}\cdot
\det_{L^{\bot }}\varepsilon ^{2}~\Rightarrow ~\det_{L}\varepsilon ^{2}=\det_{%
\mathbb{R}^{8}}\varepsilon ^{2}\det_{L^{\bot }}\varepsilon ^{2}=(-1)^{2\cdot
(8+1)}\det_{L^{\bot }}\varepsilon ^{2}\text{.}
\end{equation*}%
Let $e_{1},..,e_{8}$ be the standard base of $\mathbb{R}^{8}$. Then one base
for $L^{\bot }$ is%
\begin{equation*}
f_{1}=e_{1}+e_{2}\text{, }f_{2}=e_{3}+e_{4}\text{, }f_{3}=e_{5}+e_{6}\text{,
}f_{4}=e_{1}+e_{2}+..+e_{7}+e_{8}\text{, }
\end{equation*}%
and $\varepsilon ^{2}$ acts on it by%
\begin{equation*}
\varepsilon ^{2}\cdot f_{1}=f_{2}\text{, }\varepsilon ^{2}\cdot f_{2}=f_{3}%
\text{, }\varepsilon ^{2}\cdot f_{3}=f_{4}-f_{1}-f_{2}-f_{3}\text{, }%
\varepsilon ^{2}\cdot f_{4}=f_{4}\text{.}
\end{equation*}%
Therefore,%
\begin{equation*}
\det_{L^{\bot }}\varepsilon ^{2}=\det \left[
\begin{array}{cccc}
{\small 0} & {\small 1} & {\small 0} & {\small 0} \\
{\small 0} & {\small 0} & {\small 1} & {\small 0} \\
{\small -1} & {\small -1} & {\small -1} & {\small 1} \\
0 & 0 & 0 & {\small 1}%
\end{array}%
\right] =-1~\Rightarrow \epsilon =\det_{L}\varepsilon ^{2}=-\det_{\mathbb{R}%
^{8}}\varepsilon ^{2}=-1~\Rightarrow l=-(\varepsilon ^{2}\cdot l)
\end{equation*}%
Now, we calculate coinvariants $H_{2}(W_{n}\backslash \cup \mathcal{A};%
\mathbb{Z})_{\mathbb{Z}_{8}}$ by using the modified action $g\ast l=\det
(g)~g\cdot l$, $g\in \mathbb{Z}_{8}$ on $\mathbb{Z}_{8}$-module $H_{4}(\cup
\widehat{\mathcal{A}};\mathbb{Z})$. Let us observe that
\begin{equation*}
l\sim \varepsilon \ast l=\det (\varepsilon )\varepsilon \cdot l=-\varepsilon
\cdot l\text{ and }l\sim \varepsilon ^{2}\ast l=\det (\varepsilon
^{2})\varepsilon ^{2}\cdot l=-\det (\varepsilon ^{2})\det (\varepsilon
^{2})l=-l\text{.}
\end{equation*}%
Thus, from this relations we can conclude that $H_{2}(W_{n}\backslash \cup
\mathcal{A};\mathbb{Z})_{\mathbb{Z}_{8}}\cong \mathbb{Z}_{2}$.
\end{example}

\begin{example}
\label{Ex:Coinvariants2}Let $\mathcal{B}$ be the minimal $\mathbb{Z}_{4}$
arrangement in $\mathbb{R}^{8}=\mathbb{R}^{4}\oplus \mathbb{R}^{4}$
containing subspace
\begin{equation*}
L=\{\mathbf{x}\in \mathbb{R}%
^{8}~|~x_{1}=x_{5}=x_{1}+..+x_{4}=x_{5}+..+x_{8}=x_{3}+x_{7}=0\}\subset
W_{4}\oplus W_{4}.
\end{equation*}%
Let the action of $\mathbb{Z}_{4}=\langle \varepsilon \rangle $ be defined
by $\varepsilon \cdot
(x_{1},..,x_{8})=(x_{4},x_{1},x_{2},x_{3};x_{8},x_{5},x_{6},x_{7})$. Then $%
\det (\varepsilon )=(-1)^{4+1}(-1)^{4+1}=1$. The arrangement $\mathcal{A}$
has $4$ maximal elements $L$, $\varepsilon L$, $\varepsilon ^{2}L$, $%
\varepsilon ^{3}L$ and the Hasse diagram of the intersection poset
$P$ is like in the figure \ref{fig:Fig3}.


\begin{figure}[htb]
\centering
\includegraphics[scale=0.6]{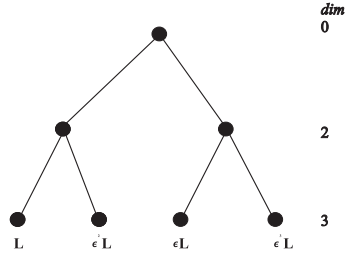}
\caption{The Hasse diagram of $\mathcal{B}$.} \label{fig:Fig3}
\end{figure}

\begin{equation*}
H_{3}(\cup \widehat{\mathcal{B}};\mathbb{Z})\cong \underset{d=0}{\overset{3}{%
\bigoplus }}\underset{V\in Q(\alpha ):\dim V=d}{\bigoplus
}H_{2n-5}(\Delta
(P_{<V})\ast \hat{V})\cong 0\oplus (\mathbb{Z}\oplus \mathbb{Z})\oplus (%
\mathbb{Z}\oplus \mathbb{Z\oplus Z}\oplus \mathbb{Z}).
\end{equation*}%
Like in the previous example the homology $H_{3}(\cup \widehat{\mathcal{B}};%
\mathbb{Z})$ has no torsion and so $H_{2}(W_{4}\oplus W_{4}\backslash \cup
\mathcal{B};\mathbb{Z})\cong H_{3}(\cup \widehat{\mathcal{B}};\mathbb{Z})$.

\noindent Let $l$, $\varepsilon \cdot l$, $\varepsilon ^{2}\cdot l$, $%
\varepsilon ^{3}\cdot l$ and $k$, $\varepsilon \cdot k$ be the generators of
$H_{3}(\cup \widehat{\mathcal{B}};\mathbb{Z})$ corresponding to spheres $%
\widehat{L}$, $\widehat{\varepsilon L}$, $\widehat{\varepsilon ^{2}L}$, $%
\widehat{\varepsilon ^{3}L}$ and $\widehat{L\cap \varepsilon ^{2}L}\ast S^{0}
$, $\widehat{\varepsilon L\cap \varepsilon ^{3}L}\ast S^{0}$ respectively.
The equality $L\cap \varepsilon ^{2}L=\varepsilon ^{2}(L\cap \varepsilon
^{2}L)$ in $\mathbb{R}^{8}$ implie the equality $k=\epsilon (\varepsilon
^{2}k)$, where $\epsilon \in \{1,-1\}$. The sign $\epsilon $ depends on the
nature of $\varepsilon ^{2}$. If $\varepsilon ^{2}$ changes the orientation
of the sphere $\widehat{L\cap \varepsilon ^{2}L}\ast S^{0}$, then $\epsilon
=-1$, otherwise $\epsilon =1$. Again, we use the same decomposition $\mathbb{%
R}^{8}=(L\cap \varepsilon ^{2}L)\oplus (L\cap \varepsilon ^{2}L)^{\bot }$
and calculate%
\begin{equation*}
\det_{\mathbb{R}^{8}}\varepsilon ^{2}=\det_{(L\cap \varepsilon
^{2}L)}\varepsilon ^{2}\cdot \det_{(L\cap \varepsilon ^{2}L)^{\bot
}}\varepsilon ^{2}~\Rightarrow ~\det_{(L\cap \varepsilon ^{2}L)}\varepsilon
^{2}=\det_{\mathbb{R}^{8}}\varepsilon ^{2}\det_{(L\cap \varepsilon
^{2}L)^{\bot }}\varepsilon ^{2}=\det_{(L\cap \varepsilon ^{2}L)^{\bot
}}\varepsilon ^{2}\text{.}
\end{equation*}%
If $e_{1},..,e_{8}$ is the standard base of $\mathbb{R}^{8}$, then one base
for $(L\cap \varepsilon ^{2}L)^{\bot }$ is%
\begin{eqnarray*}
f_{1} &=&e_{1}\text{, }f_{2}=e_{1}+e_{2}+e_{3}+e_{4}\text{, }f_{3}=e_{5}%
\text{, } \\
f_{4} &=&e_{5}+e_{6}+e_{7}+e_{8}\text{, }f_{5}=e_{3}\text{, }f_{6}=e_{7}%
\text{, }
\end{eqnarray*}%
and $\varepsilon ^{2}$ acts on it by%
\begin{equation*}
\varepsilon ^{2}\cdot f_{1}=f_{5}\text{, }\varepsilon ^{2}\cdot f_{2}=f_{2}%
\text{, }\varepsilon ^{2}\cdot f_{3}=f_{6}\text{, }\varepsilon ^{2}\cdot
f_{4}=f_{4}\text{, }\varepsilon ^{2}\cdot f_{5}=f_{1}\text{, }\varepsilon
^{2}\cdot f_{6}=f_{3}.
\end{equation*}%
Since the sign of the permutation $\left(
\begin{array}{llllll}
{\small 1} & {\small 2} & {\small 3} & {\small 4} & {\small 5} & {\small 6}
\\
{\small 5} & {\small 2} & {\small 6} & {\small 4} & {\small 1} & {\small 3}%
\end{array}%
\right) $ is one then $\det_{(L\cap \varepsilon ^{2}L)^{\bot }}\varepsilon
^{2}=1$.

\noindent In order to verify that the element $\varepsilon ^{2}$ does not
change the orientation of the sphere $\widehat{L\cap \varepsilon ^{2}L}\ast
S^{0}$, we analyze the action of $\varepsilon ^{2}$ on the homology $H_{3}(%
\widehat{L\cap \varepsilon ^{2}L}\ast S^{0};\mathbb{Z})$. Since $\widehat{%
L\cap \varepsilon ^{2}L}$ is a $2$-sphere and, we saw, $\varepsilon ^{2}$
acts trivially on it, the isomorphism%
\begin{equation*}
H_{3}(\widehat{L\cap \varepsilon ^{2}L}\ast S^{0};\mathbb{Z})\cong \tilde{H}%
_{0}(S^{0};\mathbb{Z})\otimes H_{2}(\widehat{L\cap \varepsilon ^{2}L};%
\mathbb{Z})
\end{equation*}%
instructs us that it remains to look at the action of the $\varepsilon ^{2}$
on $\tilde{H}_{0}(S^{0};\mathbb{Z})$. Keeping in mind that $S^{0}$ is the
order complex of the lower cone of the element $L\cap \varepsilon ^{2}L$, it
is obvious that $\varepsilon ^{2}$ acts on $H_{0}(S^{0};\mathbb{Z})\cong
\mathbb{Z\oplus Z}$ by permuting generators of the two copies of $\mathbb{Z}$
and remembering the orientation, i. e.%
\begin{equation*}
(x,y)~\mapsto ~\det \varepsilon ^{2}(y,x)
\end{equation*}%
Therefore, the definition of the augmentation implies that $\varepsilon ^{2}$
acts by multiplication on the reduced homology $\tilde{H}_{0}(S^{0};\mathbb{Z%
})\cong \mathbb{Z}$. Thus, $k=(\det \varepsilon ^{2})~\varepsilon ^{2}\cdot k
$.

\noindent Again we compute the coinvariants $H_{2}(W_{4}\oplus
W_{4}\backslash \cup \mathcal{B};\mathbb{Z})_{\mathbb{Z}_{4}}$ using the
modified action. The following relations
\begin{equation*}
l\sim \varepsilon ^{i}\ast l=\det (\varepsilon ^{i})\varepsilon ^{i}\cdot
l=\varepsilon ^{i}\cdot l\text{, }k\sim \varepsilon \ast k=\det (\varepsilon
)\varepsilon \cdot k=\varepsilon k\text{, }k\sim \varepsilon ^{2}\ast
k=(\det \varepsilon ^{2})\varepsilon ^{2}\cdot k=k
\end{equation*}%
imply that actually nothing "torsion-like" happens and that $%
H_{2}(W_{4}\oplus W_{4}\backslash \cup \mathcal{B};\mathbb{Z})_{\mathbb{Z}%
_{4}}\cong \mathbb{Z\oplus Z}$.
\end{example}

\subsection{\label{sec:HowToComputeTheOstructionCocycle}How to compute the
obstruction cocycle}

\noindent Finally, we are ready to give an algorithm for computing the
cohomology class of the obstruction cocycle of the map in general position.
Let us assume that additional hyper arrangement $\mathcal{J}$ is already
chosen.

\textit{Step 1:} Let $S^{3}$ be $\mathbb{Q}_{4n}$ the simplicial complex
earlier defined. Define $\mathbb{Q}_{4n}$-map $h\,:S^{3}\rightarrow W_{n}$
by defining them on $0$-skeleton. It is enough to define an image of a
single vertex, because everything else is defined by the equivariant request.

\textit{Step 2:}\textbf{\ }Find all simplexes $\sigma
_{ij}=[a_{i},a_{i+1};b_{j},b_{j+1}]$ such that
\begin{equation*}
h([a_{i},a_{i+1};b_{j},b_{j+1}])\cap (L(\alpha )\cap \mathcal{J})=\{pt.\}%
\text{.}
\end{equation*}%
The intersection can't have more than one point, because then at least the
whole interval will be in the intersection. This would mean that the map $h$
is not in the general position.

\textit{Step 3:} With the help of the cellular map between two introduced $%
\mathbb{Q}_{4n}$ cellular structures on $S^{3}$, and previous step, count
the following sets%
\begin{equation*}
h(e)\cap (\cup \mathcal{A}(\mathcal{J},\alpha ))\text{ and }h^{-1}(h(e)\cap
(\cup \mathcal{A}(\mathcal{J},\alpha )))\subset e\text{.}
\end{equation*}

\textit{Step 4:} Let us assume that every element $y\in h(e)\cap (\cup
\mathcal{A}(\mathcal{J},\alpha ))$ is contained in just one and only one
maximal element $L_{y}$ of the arrangement $\mathcal{A}$. Then
\begin{equation}
c_{\mathbb{Q}_{4n}}(h)(e)=\sum_{x\in h^{-1}(h(e)\cap (\cup \mathcal{A}(%
\mathcal{J},\alpha )))}\mathrm{I}(e,L_{h(x)})\left\Vert h(x)\right\Vert
\label{eq:ObstructionCocycle}
\end{equation}%
where $\mathrm{I}(h(\theta ),L_{h(x)})$ is the intersection number of the
oriented cell $e$ and appropriate oriented element $L_{y}$ ($L_{y}\cap
h(\theta )=\{y\}$) of the arrangement $\mathcal{A}(J,\alpha )$.

\noindent If some $y\in h(e)\cap (\cup \mathcal{A}(\mathcal{J},\alpha ))$
belongs to codimension one intersection of maximal elements $%
L_{y}^{(1)},..,L_{y}^{(k)}$ the formulas is unchanged except the class $%
\left\Vert h(x)\right\Vert $ is a broken point class. The real trouble with
this situation is that we need an extra effort to identify this broken point
class. As we mentioned in section \ref{sec:pointClasses} the broken point
class depends on the embedding of the simplex (or its linear span) which
intersects arrangement in $y\in h(e)\cap (\cup \mathcal{A}(\mathcal{J}%
,\alpha ))$.

\textit{Step 5:} Compute $H_{n-4}(\cup \widehat{\mathcal{A}}(J,\alpha );%
\mathbb{Z})$ using decomposition (\ref{Z-Z homology decomposition - 1}). If
there are no torsion we actually calculated $H_{2}(W_{n}\backslash \cup
\mathcal{A}(\mathcal{J},\alpha );\mathbb{Z})$.

\textit{Step 6:} Compute the coinvariants $H_{2}(W_{n}\backslash \cup
\mathcal{A}(\mathcal{J},\alpha );\mathbb{Z})_{\mathbb{Q}_{4n}}$ by working
in $\mathbb{Q}_{4n}$-module $H_{n-4}(\cup \widehat{\mathcal{A}}(\mathcal{J}%
,\alpha );\mathbb{Z})$ using modified action - isomorphism (\ref{isomorphism
of the modified action}).

\textit{Step 7:} Express the obstruction element $c_{\mathbb{Q}%
_{4n}}(h)(e)\in H_{2}(W_{n}\backslash \cup \mathcal{A}(\mathcal{J},\alpha );%
\mathbb{Z})$ as the element of the group $\mathrm{Hom}(H_{n-4}(\cup \widehat{%
\mathcal{A}}(\mathcal{J},\alpha );\mathbb{Z}),\mathbb{Z})$ via the
isomorphism (\ref{Isomorphisms Complement - Arrangement}).

\textit{Step 8:} Identify the class of the obstruction cocycle $c_{\mathbb{Q}%
_{4n}}(h)(e)$ in the group
\begin{equation*}
H_{n-4}(\cup \widehat{\mathcal{A}}(\mathcal{J},\alpha ),\mathbb{Z})/\sim
\text{.}
\end{equation*}%
If it is not zero we proved that the $\mathbb{Q}_{4n}$-map in question can
not exist. Thus, propositions \ref{prop:VezaProblem-Ekvivarijantan} and \ref%
{prop:PrelazakNaNovuGrupu} imply that the appropriate fan partition exists.
If it is zero, then the whole effort was in vane.

\section{\textsf{Computations and proof of theorem \protect\ref{th:main1}}}

\noindent The proof of the theorem \ref{th:main1} has two stages

\begin{itemize}
\item The proof for the wide class of special cases; in particular we prove
theorem for all $\alpha =(\tfrac{a}{n},\tfrac{a+b}{n},\tfrac{b}{n})\in \frac{%
1}{n}\,\mathbb{N}^{3}\subseteq \mathbb{Q}^{3}$ such that $2a+2b=n$, $a,b\geq
1$. This proof goes along the lines described in section \ref%
{sec:HowToComputeTheOstructionCocycle}.

\item The limit argument extend the result from the class of special cases
to the whole class of triples $\alpha =(a,a+b,b)\in \mathbb{R}_{>0}^{3}$, $%
2a+2b=1$.
\end{itemize}

\subsection{Finding hyper arrangement $\mathcal{J}$}

\noindent Before we try to prove the main theorem in steps described in the
section \ref{sec:HowToComputeTheOstructionCocycle}, let us make the
fundamental step by defining the additional hyper arrangement $\mathcal{J}$
in $\mathbb{R}^{n}$. Let us define hyperplanes $H_{1}$, $H_{2}$, $K$ and
half-spaces $K^{+}$, $K^{-}$, by%
\begin{equation*}
\begin{array}{l}
H_{1}=\{\mathbf{x}\in \mathbb{R}^{n}|\text{ }%
(a+b)(x_{a}-x_{2a+b}+x_{1}-x_{a+b+1})+x_{a+1}-x_{2a+b+1}+x_{n}-x_{a+b}=0\},
\\
H_{2}=\{\mathbf{x}\in \mathbb{R}^{n}|\text{ }x_{a+1}+..+x_{a+b}=0\},K=\{%
\mathbf{x}\in \mathbb{R}^{n}|\text{~}x_{1}+..+x_{a+b}=0\} \\
K^{+}=\{\mathbf{x}\in \mathbb{R}^{n}|\text{ }x_{1}+..+x_{a+b}\geq 0\}\text{
and }K^{-}=\{\mathbf{x}\in \mathbb{R}^{n}|\text{~}x_{1}+..+x_{a+b}\leq 0\}.%
\end{array}%
\end{equation*}%
The hyper arrangement $\mathcal{J}$ we would like to consider is%
\begin{equation*}
\mathcal{J}=\{H_{1}\cap K_{1}^{+},~H_{2}\}.
\end{equation*}%
Let us consider some properties of this arrangement which will produce
crucial arguments in the proof of the main theorem. The first property is
that $\varepsilon ^{a+b}H_{1}=H_{1}$ and $\varepsilon ^{a+b}$ does change
the orientation of the subspace $H_{1}^{\bot }$. Indeed, let $%
e=(a+b)(e_{a}-e_{2a+b}+e_{1}-e_{a+b+1})+e_{a+1}-e_{2a+b+1}+e_{n}-e_{a+b}$ be
a base vector of $H_{1}^{\bot }$, then $\varepsilon ^{a+b}\cdot
e=(a+b)(e_{2a+b}-e_{a}+e_{a+b+1}-e_{1})+e_{2a+b+1}-e_{a+1}+e_{a+b}-e_{n}$.
The other property is that $\varepsilon ^{a+b}(K^{+}\cap W_{n})=K^{-}\cap
W_{n}$.

\noindent Thus the arrangement $\mathcal{A}(\mathcal{J},\alpha )$ is the
minimal $\mathbb{Q}_{4n}$-arrangement containing half-subspace $L_{1}^{\ast
}=L(\alpha )\cap H_{1}\cap K^{+}$ and subspace $L_{2}^{\ast }=L(\alpha )\cap
H_{2}$ defined by
\begin{eqnarray*}
L_{1}^{\ast } &:&\left\{
\begin{tabular}{l}
$x_{1}+..+x_{a}=x_{a+1}+..+x_{2a+b}=x_{2a+b+1}+..+x_{n}=0,~x_{1}+..+x_{a+b}%
\geq 0,$ \\
$(a+b)(x_{a}-x_{2a+b}+x_{1}-x_{a+b+1})+x_{a+1}-x_{2a+b+1}+x_{n}-x_{a+b}=0.$%
\end{tabular}%
\right. \\
L_{2}^{\ast }
&:&x_{1}+x_{2}+..+x_{a}=x_{a+1}+..+x_{a+b}=x_{a+b+1}+..+x_{2a+b}=x_{2a+b+1}+..+x_{n}=0.
\end{eqnarray*}%
The introduced setting provides us with the following facts:

\textit{(i)} The intersection
\begin{equation*}
I=L_{1}^{\ast }\cap \varepsilon ^{a+b}L_{1}^{\ast }\cap \varepsilon
^{a}jL_{1}^{\ast }\cap \varepsilon ^{2a+b}jL_{1}^{\ast }\cap L_{2}^{\ast }
\end{equation*}%
is a linear subspace of codimension one in both $L_{1}^{\ast }=L(\alpha
)\cap H_{1}\cap K^{+}$ and $L_{2}^{\ast }=L(\alpha )\cap H_{2}$;

\textit{(ii)} The following set of equalities stands%
\begin{equation*}
\varepsilon ^{a+b}(L_{1}^{\ast }\cap \varepsilon ^{a+b}L_{1}^{\ast
})=L_{1}^{\ast }\cap \varepsilon ^{a+b}L_{1}^{\ast },~\varepsilon
^{a+b}(\varepsilon ^{a}jL_{1}^{\ast }\cap \varepsilon ^{2a+b}jL_{1}^{\ast
})=\varepsilon ^{a}jL_{1}^{\ast }\cap \varepsilon ^{2a+b}jL_{1}^{\ast
},~\varepsilon ^{a+b}L_{2}^{\ast }=L_{2}^{\ast }.
\end{equation*}

\textit{(iii)} The element $\varepsilon ^{a+b}$ changes the orientation on $%
I^{\bot }$.

\begin{figure}[htb]
\centering
\includegraphics[scale=1]{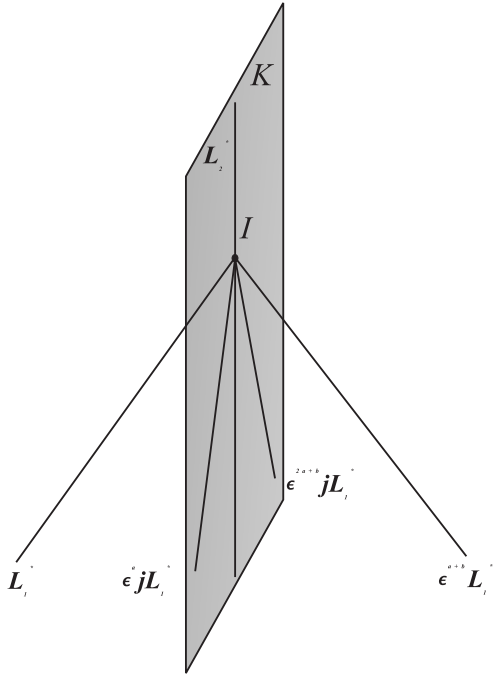}
\caption{The part of the arrangement
$\mathcal{A}(\mathcal{J},\alpha )$} \label{fig:Fig6}
\end{figure}

\noindent These properties will be the key arguments that the appropriate
group of coinvariants $H_{2}(W_{n}\backslash \cup \mathcal{A}(\mathcal{J}%
,\alpha ),\mathbb{Z})_{\mathbb{Q}_{4n}}$ have the nontrivial
torsion part. The part of our arrangement $\cup
\mathcal{A}(\mathcal{J},\alpha )$ can be pictured like in the
figure \ref{fig:Fig6}.


\subsection{Step 1}

\noindent\ Let us recall that we substituted the sphere $S^{3}$ with the
simplicial complex $P_{2n}^{(1)}\ast P_{2n}^{(2)}$, where $P_{2n}$ is the
regular $2n$-gon. Earlier we have denoted vertices of two copies of $P_{2n}$
by $a_{1},..,a_{2n}$ and $b_{1},..,b_{2n}$, respectively. To simplify
notation in calculation ahead let $t=a_{1}$. Let $e_{1},..,e_{n}$ be the
standard base of $\mathbb{R}^{n}$ and $u_{i}=e_{i}-\tfrac{1}{n}%
\sum\limits_{j=1}^{n}e_{j}$, $i\in \{1,..,n\}$.

\noindent We define the $\mathbb{Q}_{4n}$-map $h:S^{3}\rightarrow
W_{n}\subset \mathbb{R}^{n}$ which is in the general position by defining it
on the vertex $t$ by $h(t)=u_{1}$. Then the request that $h$ is $\mathbb{Q}%
_{4n}$-map imply
\begin{equation*}
\begin{array}{l}
h(\varepsilon ^{i}t)=\varepsilon ^{i}h(t)=\varepsilon ^{i\
\textrm{mod}\ %
n}h(t)=u_{i\ \textrm{mod}\ n+1} \\
h(jt)=jh(t)=u_{n} \\
h(\varepsilon ^{i}jt)=\varepsilon ^{i}jh(t)=\varepsilon ^{i\
\textrm{mod}\ %
n}jh(t)=u_{i\ \textrm{mod}\ n}.%
\end{array}%
\end{equation*}%
In the future we will omit the "$\textrm{mod}\ n$" part in the
indexes on the
right hand side, i.e. all the indexes in $W_{n}$ are calculated $\textrm{mod}\ n$%
.

\subsection{Step 2\textbf{\ }}

\noindent It is not hard to see from the definition that the image
of the map $h$ is
$h(S^{3})=\bigcup\limits_{i,j=1}^{n}[u_{i},u_{i+1}]\ast \lbrack
u_{j},u_{j+1}]\subset W_{n}$, remembering that $n+1$ is actually
$1$. Let us list all the simplexes of the form
$[u_{i},u_{i+1}]\ast \lbrack u_{j},u_{j+1}]\equiv \lbrack
u_{i},u_{i+1};u_{j},u_{j+1}]$ which intersect
the subspace%
\begin{equation*}
L(\alpha )=\{\mathbf{x}\in W_{n}|~x_{1}+\ldots +x_{a}=0,x_{a+1}+\ldots
+x_{2a+b}=0,x_{2a+b+1}+\ldots +x_{n}=0\}\text{.}
\end{equation*}%
With a little linear algebra and combinatorics we see that the only
simplexes from the image $h(S^{3})$ that intersects $L(\alpha )$ are:%
\begin{equation*}
\begin{tabular}{|l|l|l|}
\hline
\ \ {\small Simplex} & {\small Intersection with }${\small L(\alpha )}$ &
{\small For }${\small r\in }$ \\ \hline
{\small 1.}${\small [u}_{a}{\small ,u}_{a+1}{\small ;u}_{2a+b}{\small ,u}%
_{2a+b+1}{\small ]}$ & ${\small \{}\tfrac{a}{n}{\small u}_{a}{\small +}%
\tfrac{\beta }{n}{\small u}_{a+1}{\small +}\tfrac{\gamma }{n}{\small u}%
_{2a+b}{\small +}\tfrac{b}{n}{\small u}_{2a+b+1}{\small \}}$ &  \\ \hline
{\small 2.}${\small [u}_{a}{\small ,u}_{a+1}{\small ;u}_{r}{\small ,u}_{r+1}%
{\small ]}$ & $\{\tfrac{a}{n}{\small u}_{a}{\small +}\tfrac{a+b}{n}{\small u}%
_{a+1}{\small +}\tfrac{\gamma }{n}{\small u}_{r}{\small +}\tfrac{\delta }{n}%
{\small u}_{r+1}\}$ & ${\small [2a+b+1,n-1]}$ \\ \hline
{\small 3.}${\small [u}_{r}{\small ,u}_{r+1}{\small ;u}_{2a+b}{\small ,u}%
_{2a+b+1}{\small ]}$ & $\{\tfrac{\alpha }{n}{\small u}_{r}{\small +}\tfrac{%
\beta }{n}{\small u}_{r+1}{\small +}\tfrac{a+b}{n}{\small u}_{2a+b}{\small +}%
\tfrac{b}{n}{\small u}_{2a+b+1}{\small \}}$ & ${\small [1,a-1]}$ \\ \hline
{\small 4.}${\small [u}_{a}{\small ,u}_{a+1}{\small ;u}_{n}{\small ,u}_{1}%
{\small ]}$ & $\{\tfrac{\alpha }{n}{\small u}_{a}{\small +}\tfrac{a+b}{n}%
{\small u}_{a+1}{\small +}\tfrac{b}{n}{\small u}_{n}{\small +}\tfrac{\delta
}{n}{\small u}_{1}{\small \}}$ &  \\ \hline
{\small 5.}${\small [u}_{2a+b}{\small ,u}_{2a+b+1}{\small ;u}_{n}{\small ,u}%
_{1}{\small ]}$ & $\{\tfrac{a+b}{n}{\small u}_{2a+b}{\small +}\tfrac{\beta }{%
n}{\small u}_{2a+b+1}{\small +}\tfrac{\gamma }{n}{\small u}_{n}{\small +}%
\tfrac{a}{n}{\small u}_{1}{\small \}}$ &  \\ \hline
{\small 6.}${\small [u}_{r}{\small ,u}_{r+1}{\small ;u}_{n}{\small ,u}_{1}%
{\small ]}$ & $\{\tfrac{\alpha }{n}{\small u}_{r}{\small +}\tfrac{\beta }{n}%
{\small u}_{r+1}{\small +}\tfrac{b}{n}{\small u}_{n}{\small +}\tfrac{a}{n}%
{\small u}_{1}{\small \}}$ & ${\small [a+1,2a+b-1]}$ \\ \hline
\end{tabular}%
\end{equation*}%
Whit the assumption that $a,b\geq 1$ we analize two cases.

\textit{(1)} If we introduce the hyperplane $H_{1}$, then the only simplex
that intersects $L(\alpha )\cap H_{1}$ are%
\begin{equation*}
\rho _{1}={\small [u}_{a}{\small ,u}_{a+1}{\small ;u}_{2a+b}{\small ,u}%
_{2a+b+1}{\small ]}\text{, }\rho _{2}={\small [u}_{a+b}{\small ,u}_{a+b+1}%
{\small ;u}_{n}{\small ,u}_{1}{\small ]}\text{, }\rho _{3}={\small [u}%
_{a+b+1}{\small ,u}_{a+b+2}{\small ;u}_{n}{\small ,u}_{1}{\small ]}\text{.}
\end{equation*}%
Precisely,%
\begin{eqnarray*}
(L(\alpha )\cap H_{1})\cap \rho _{1} &=&\{\tfrac{a}{n}u_{a}+\tfrac{b}{n}%
u_{a+1}+\tfrac{a}{n}u_{2a+b}+\tfrac{b}{n}u_{2a+b+1}\} \\
(L(\alpha )\cap H_{1})\cap \rho _{2} &=&\{\tfrac{b}{n}u_{a+b}+\tfrac{a}{n}%
u_{a+b+1}+\tfrac{b}{n}u_{n}+\tfrac{a}{n}u_{1}\} \\
(L(\alpha )\cap H_{1})\cap \rho _{3} &=&\{\tfrac{a+\frac{b}{a+b}}{n}{\small u%
}_{a+b+1}{\small +}\tfrac{b-\frac{b}{a+b}}{n}{\small u}_{a+b+2}{\small +}%
\tfrac{b}{n}{\small u}_{n}{\small +}\tfrac{a}{n}{\small u}_{1}\}.
\end{eqnarray*}%
But when we add the inequality condition $x_{1}+..+x_{a+b}\geq 0$, the
intersection of $L_{1}^{\ast }=L(\alpha )\cap H_{1}\cap K^{+}$ with the
simplex $\rho _{3}$ vanishes. Thus there are only two simplexes
\begin{equation*}
{\small [u}_{a}{\small ,u}_{a+1}{\small ;u}_{2a+b}{\small ,u}_{2a+b+1}%
{\small ]}\text{ and }{\small [u}_{a+b}{\small ,u}_{a+b+1}{\small ;u}_{n}%
{\small ,u}_{1}{\small ]}
\end{equation*}%
which intersects $L_{1}^{\ast }$, each of them in the single interior point,%
\begin{equation*}
v=\tfrac{a}{n}u_{a}+\tfrac{b}{n}u_{a+1}+\tfrac{a}{n}u_{2a+b}+\tfrac{b}{n}%
u_{2a+b+1}\in L_{1}^{\ast }\text{ and }w=\tfrac{b}{n}u_{a+b}+\tfrac{a}{n}%
u_{a+b+1}+\tfrac{b}{n}u_{n}+\tfrac{a}{n}u_{1}\in L_{1}^{\ast }\text{.}
\end{equation*}%
Thus, $h(S^{3})\cap L_{1}^{\ast }=\{v,w\}$.

\textit{(2)} If we introduce the hyperplane $H_{2}$ instead the hyperplane $%
H_{1}$, there are only two simplexes from $h(S^{3})$ that intersects $%
L_{2}^{\ast }=L(\alpha )\cap H_{2}$. Those are
\begin{equation*}
\rho _{1}={\small [u}_{a}{\small ,u}_{a+1}{\small ;u}_{2a+b}{\small ,u}%
_{2a+b+1}{\small ]}\text{ and }\rho _{2}={\small [u}_{a+b}{\small ,u}_{a+b+1}%
{\small ;u}_{n}{\small ,u}_{1}{\small ]}
\end{equation*}%
and by another miracle they intersect $L_{2}^{\ast }$ in the same points as
in the previous case, $h(S^{3})\cap L_{2}^{\ast }=\{v,w\}$.

\begin{conclusion}
This means that in order to analize the set $h^{-1}(h(e)\cap (\cup \mathcal{A%
}(\mathcal{J},\alpha )))\subset e$ we only have to track down the pre images
of $v$ and $w$.
\end{conclusion}

\subsection{Step 3}

\noindent There are 16 simplexes in the sphere $P_{2n}^{(1)}\ast P_{2n}^{(2)}
$ which belong to the $h$ inverse image of the simplexes ${\small [u}_{a}%
{\small ,u}_{a+1}{\small ;u}_{2a+b}{\small ,u}_{2a+b+1}{\small ]}$ and $%
{\small [u}_{a+b}{\small ,u}_{a+b+1}{\small ;u}_{n}{\small ,u}_{1}{\small ]}$%
. These simplexes are
\begin{equation*}
\begin{array}{ll}
\sigma _{1}=[\epsilon ^{a-1}t,\epsilon ^{a}t;\epsilon ^{2a+b}jt,\epsilon
^{2a+b+1}jt], & \sigma _{2}=[\epsilon ^{n+a-1}t,\epsilon ^{n+a}t;\epsilon
^{2a+b}jt,\epsilon ^{2a+b+1}jt], \\
\sigma _{3}=[\epsilon ^{a-1}t,\epsilon ^{a}t;\epsilon ^{n+2a+b}jt,\epsilon
^{n+2a+b+1}jt], & \sigma _{4}=[\epsilon ^{n+a-1}t,\epsilon ^{n+a}t;\epsilon
^{n+2a+b}jt,\epsilon ^{n+2a+b+1}jt], \\
\sigma _{5}=[\epsilon ^{2a+b-1}t,\epsilon ^{2a+b}t;\epsilon ^{a}jt,\epsilon
^{a+1}jt], & \sigma _{6}=[\epsilon ^{n+2a+b-1}t,\epsilon ^{n+2a+b}t;\epsilon
^{a}jt,\epsilon ^{a+1}jt], \\
\sigma _{7}=[\epsilon ^{2a+b-1}t,\epsilon ^{2a+b}t;\epsilon
^{n+a}jt,\epsilon ^{n+a+1}jt], & \sigma _{8}=[\epsilon ^{n+2a+b-1}t,\epsilon
^{n+2a+b}t;\epsilon ^{n+a}jt,\epsilon ^{n+a+1}jt], \\
\theta _{1}=[\epsilon ^{a+b-1}t,\epsilon ^{a+b}t;jt,\epsilon jt], & \theta
_{2}=[\epsilon ^{n+a+b-1}t,\epsilon ^{n+a+b}t;jt,\epsilon jt], \\
\theta _{3}=[\epsilon ^{a+b-1}t,\epsilon ^{a+b}t;\epsilon ^{n}jt,\epsilon
^{n+1}jt], & \theta _{4}=[\epsilon ^{n+a+b-1}t,\epsilon ^{n+a+b}t;\epsilon
^{n}jt,\epsilon ^{n+1}jt], \\
\theta _{5}=[\epsilon ^{2n-1}t,t;\epsilon ^{a+b}jt,\epsilon ^{a+b+1}jt], &
\theta _{6}=[\epsilon ^{2n-1}t,t;\epsilon ^{n+a+b}jt,\epsilon ^{n+a+b+1}jt],
\\
\theta _{7}=[\epsilon ^{n-1}t,\epsilon ^{n}t;\epsilon ^{a+b}jt,\epsilon
^{a+b+1}jt], & \theta _{8}=[\epsilon ^{n-1}t,\epsilon ^{n}t;\epsilon
^{n+a+b}jt,\epsilon ^{n+a+b+1}jt]%
\end{array}%
\end{equation*}%
By some miracle all these simplexes are in the orbit of the single simplex $%
\sigma $ in the maximal cell $e$. Moreover, the pre-images of intersection
points with $L^{\ast }$%
\begin{equation*}
v=\tfrac{a}{n}u_{a}+\tfrac{b}{n}u_{a+1}+\tfrac{a}{n}u_{2a+b}+\tfrac{b}{n}%
u_{2a+b+1}\text{ and }w=\tfrac{b}{n}u_{a+b}+\tfrac{a}{n}u_{a+b+1}+\tfrac{b}{n%
}u_{n}+\tfrac{a}{n}u_{1}
\end{equation*}%
are all in the orbit of two points
\begin{equation*}
v^{\ast }=\tfrac{a}{n}\epsilon ^{a+b-1}t+\tfrac{b}{n}\epsilon ^{a+b}t+\tfrac{%
a}{n}jt+\tfrac{b}{n}\epsilon jt\text{ and }w^{\ast }=\tfrac{b}{n}\epsilon
^{a+b-1}t+\tfrac{a}{n}\epsilon ^{a+b}t+\tfrac{b}{n}jt+\tfrac{a}{n}\epsilon jt
\end{equation*}%
in the simplex $\sigma $. In particular,%
\begin{equation*}
\begin{array}{ll}
\sigma =[\epsilon ^{a+b-1}t,\epsilon ^{a+b}t;jt,\epsilon jt] & =\epsilon
j\epsilon ^{-a+1}\sigma _{1}=\epsilon ^{-2a-b}\sigma _{2}=\epsilon
^{-n-2a-b}\sigma _{3}=\epsilon j\epsilon ^{-n-a+1}\sigma _{4} \\
& =\epsilon ^{-a}\sigma _{5}=\epsilon j\epsilon ^{-n-2a-b+1}\sigma
_{6}=\epsilon j\epsilon ^{-2a-b+1}\sigma _{7}=\epsilon ^{-a-n}\sigma _{8} \\
& =\theta _{1}=\epsilon j\epsilon ^{-n-a-b+1}\theta _{2}=\epsilon j\epsilon
^{-a-b+1}\theta _{3}=\epsilon ^{-n}\theta _{4} \\
& =\epsilon j\epsilon \theta _{5}=\epsilon ^{-n-a-b}\theta _{6}=\epsilon
^{-a-b}\theta _{7}=\epsilon j\epsilon ^{-n+1}\theta _{8}%
\end{array}%
\end{equation*}%
To illustrate the fact that points $v^{\ast }$ and $w^{\ast }$ from the cell
$e$ are mapped in the arrangement, we track the change of the barycentric
coordinates of points $v$, $w$ and $v^{\ast }$, $w^{\ast }$ along the
actions. For example,

\textit{(A)} From $\sigma =[\epsilon jt,jt;\epsilon ^{a+b-1}t,\epsilon
^{a+b}t]=\epsilon j\epsilon ^{-a+1}\sigma _{1}=\epsilon ^{a}j\sigma _{1}$,
we have
\begin{eqnarray*}
h(\sigma _{1})\cap L_{1}^{\ast } &=&\{\tfrac{a}{n}u_{a}+\tfrac{b}{n}u_{a+1}+%
\tfrac{a}{n}u_{2a+b}+\tfrac{b}{n}u_{2a+b+1}\} \\
h(\epsilon ^{a}j\sigma _{1})\cap \epsilon ^{a}jL_{1}^{\ast } &=&\{\tfrac{a}{n%
}u_{1}+\tfrac{b}{n}u_{n}+\tfrac{a}{n}u_{a+b+1}+\tfrac{b}{n}u_{a+b}\} \\
h^{-1}\left( h(\epsilon ^{a}j\sigma _{1})\cap \epsilon ^{a}jL_{1}^{\ast
}\right) &=&\{\tfrac{a}{n}\epsilon jt+\tfrac{b}{n}jt+\tfrac{a}{n}\epsilon
^{a+b}t+\tfrac{b}{n}\epsilon ^{a+b-1}t\}=\{w^{\ast }\}.
\end{eqnarray*}

\textit{(B)} From $\sigma =[\epsilon ^{a+b-1}t,\epsilon ^{a+b}t;jt,\epsilon
jt]=\epsilon ^{-2a-b}\sigma _{2}$, we conclude that
\begin{eqnarray*}
h(\sigma _{2})\cap L_{1}^{\ast } &=&\{\tfrac{a}{n}u_{a}+\tfrac{b}{n}u_{a+1}+%
\tfrac{a}{n}u_{2a+b}+\tfrac{b}{n}u_{2a+b+1}\} \\
h(\epsilon ^{-2a-b}\sigma _{2})\cap \epsilon ^{-2a-b}L_{1}^{\ast } &=&\{%
\tfrac{a}{n}u_{a+b}+\tfrac{b}{n}u_{a+b+1}+\tfrac{a}{n}u_{n}+\tfrac{b}{n}%
u_{a+b}\} \\
h^{-1}\left( h(\epsilon ^{-2a-b}\sigma _{2})\cap \epsilon
^{-2a-b}L_{1}^{\ast }\right) &=&\{\tfrac{a}{n}\epsilon ^{a+b-1}t+\tfrac{b}{n}%
\epsilon ^{a+b}t+\tfrac{a}{n}jt+\tfrac{b}{n}\epsilon jt\}=\{v^{\ast }\};
\end{eqnarray*}

\textit{(C)} From $\sigma =[\epsilon ^{a+b-1}t,\epsilon ^{a+b}t;jt,\epsilon
jt]=\theta _{1}$, it easily follows

\begin{eqnarray*}
h(\theta _{1})\cap L_{1}^{\ast } &=&\{\tfrac{b}{n}u_{a+b}+\tfrac{a}{n}%
u_{a+b+1}+\tfrac{b}{n}u_{n}+\tfrac{a}{n}u_{1}\} \\
h^{-1}\left( h(\theta _{1})\cap L_{1}^{\ast }\right) &=&\{\tfrac{b}{n}%
\epsilon ^{a+b-1}t+\tfrac{a}{n}\epsilon ^{a+b}t+\tfrac{b}{n}jt+\tfrac{a}{n}%
\epsilon jt\}=\{w^{\ast }\}.
\end{eqnarray*}

\begin{conclusion}
Therefore,
\begin{equation*}
\mathrm{card}(h(e)\cap (\cup \mathcal{A}(J,\alpha )))=\mathrm{card}%
(h^{-1}(h(e)\cap (\cup \mathcal{A}(J,\alpha ))))=2
\end{equation*}%
and
\begin{equation*}
h^{-1}(h(e)\cap (\cup \mathcal{A}(J,\alpha )))=\{v^{\ast },w^{\ast }\}\text{.%
}
\end{equation*}
\end{conclusion}

\subsection{Step 4}

\noindent Before analyzing the obstruction cocycle, let us observe that%
\begin{equation*}
\begin{array}{lll}
h(v^{\ast })=\epsilon ^{b}\cdot v\in  & \epsilon ^{-2a-b}(L_{1}^{\ast }\cap
L_{2}^{\ast })\cap \epsilon ^{-a}(L_{1}^{\ast }\cap L_{2}^{\ast })\cap
\epsilon j\epsilon ^{-a-b+1}(L_{1}^{\ast }\cap L_{2}^{\ast })\cap \epsilon
j\epsilon (L_{1}^{\ast }\cap L_{2}^{\ast }) & = \\
& \epsilon ^{b}(L_{1}^{\ast }\cap L_{2}^{\ast })\cap \epsilon
^{2b+a}(L_{1}^{\ast }\cap L_{2}^{\ast })\cap \epsilon ^{a+b}j(L_{1}^{\ast
}\cap L_{2}^{\ast })\cap j(L_{1}^{\ast }\cap L_{2}^{\ast }) & = \\
& \epsilon ^{b}(L_{1}^{\ast }\cap \epsilon ^{b+a}L_{1}^{\ast }\cap \epsilon
^{a}jL_{1}^{\ast }\cap \epsilon ^{2a+b}jL_{1}^{\ast }\cap L_{2}^{\ast }) &
\end{array}%
\end{equation*}%
and similarly%
\begin{equation*}
\begin{array}{lll}
h(w^{\ast })=w\in  & (L_{1}^{\ast }\cap L_{2}^{\ast })\cap \epsilon
^{-a-b}(L_{1}^{\ast }\cap L_{2}^{\ast })\cap \epsilon j\epsilon
^{-a+1}(L_{1}^{\ast }\cap L_{2}^{\ast })\cap \epsilon j\epsilon
^{-2a-b+1}(L_{1}^{\ast }\cap L_{2}^{\ast }) & = \\
& (L_{1}^{\ast }\cap L_{2}^{\ast })\cap \epsilon ^{a+b}(L_{1}^{\ast }\cap
L_{2}^{\ast })\cap \epsilon ^{a}j(L_{1}^{\ast }\cap L_{2}^{\ast })\cap
\epsilon ^{2a+b}j(L_{1}^{\ast }\cap L_{2}^{\ast }) & = \\
& L_{1}^{\ast }\cap \epsilon ^{b+a}L_{1}^{\ast }\cap \epsilon
^{a}jL_{1}^{\ast }\cap \epsilon ^{2a+b}jL_{1}^{\ast }\cap L_{2}^{\ast } &
\end{array}%
\end{equation*}%
Then the equality \ref{eq:ObstructionCocycle} applied in this situation
implies that the obstruction cocycle is
\begin{equation*}
c_{\mathbb{Q}_{4n}}(h)(e)=\left\Vert h(v^{\ast })\right\Vert +\left\Vert
h(w^{\ast })\right\Vert
\end{equation*}%
where $\tau _{1},\tau _{2}\in \{+1,-1\}$. In this situation classes $%
\left\Vert h(v^{\ast })\right\Vert $ and $\left\Vert h(w^{\ast })\right\Vert
$ are broken point classes determined by the $h$ embedding of the boundary
of the simplex $\sigma $, i.e. by the linear subspace $\mathrm{span}%
\{u_{a+b},u_{a+b+1},u_{n},u_{1}\}$.

\noindent Since, $h(v^{\ast })=\epsilon ^{b}\cdot v$ and $h(w^{\ast })=w$
the obstruction cocycle is
\begin{equation*}
c_{\mathbb{Q}_{4n}}(h)(e)=\epsilon ^{b}\cdot \left\Vert v\right\Vert
+\left\Vert w\right\Vert \text{.}
\end{equation*}%
where the first broken point class is determined by the linear span of the
simplex
\begin{equation*}
\epsilon
^{-b}[u_{a+b},u_{a+b+1},u_{n},u_{1}]=[u_{a},u_{a+1},u_{2a+b},u_{2a+b+1}].
\end{equation*}
Also, observe that the equality%
\begin{equation*}
\epsilon ^{a}j\cdot v=\epsilon ^{a}j\cdot \left( \tfrac{a}{n}u_{a}+\tfrac{b}{%
n}u_{a+1}+\tfrac{a}{n}u_{2a+b}+\tfrac{b}{n}u_{2a+b+1}\right) =\tfrac{a}{n}%
u_{1}+\tfrac{b}{n}u_{n}+\tfrac{a}{n}u_{a+b+1}+\tfrac{b}{n}u_{a+b}=w
\end{equation*}%
and the fact that the permutation $(4321)$ is even implies%
\begin{equation*}
\left\Vert w\right\Vert =\epsilon ^{a}j\cdot \left\Vert v\right\Vert \text{.}
\end{equation*}

\begin{conclusion}
The cohomology class of the obstruction cocycle, as we have seen, lives in
the group of coinvariants $H_{2}(W_{n}\backslash \cup \mathcal{A}(\mathcal{J}%
,\alpha );\mathbb{Z})_{\mathbb{Q}_{4n}}$. Thus, instead of the original
obstruction cocycle $c_{\mathbb{Q}_{4n}}(h)(e)=\epsilon ^{b}\cdot \left\Vert
v\right\Vert +\left\Vert w\right\Vert =(\epsilon ^{b}+\epsilon ^{a}j)\cdot
\left\Vert v\right\Vert $ we can analize the cocycle%
\begin{equation*}
c^{\prime }=2\left\Vert v\right\Vert
\end{equation*}%
where the broken point class $\left\Vert v\right\Vert $ is determined by the
linead span of the simplex%
\begin{equation*}
\lbrack u_{a},u_{a+1},u_{2a+b},u_{2a+b+1}].
\end{equation*}
\end{conclusion}

\subsection{Step 5}

\noindent Now we will find the complete second homology $H_{2}(W_{n}%
\backslash \cup \mathcal{A}(\mathcal{J},\alpha );\mathbb{Z})$ via the
isomorphism (\ref{hom decomposition - 1}). In order to do so we describe the
first two lower levels of the Hasse diagram of the intersection poset $%
P(\alpha )$ of the arangement $\mathcal{A}(\mathcal{J},\alpha )$. With a
little linear algebra and mentioned properties of the added arrangement $%
\mathcal{J}$, it can be seen that the $n-4$ and $n-5$ levels of the Hasse
diagram of the intersection poset $P(\alpha )$ are like in the figure \ref%
{fig:Fig7}.


\begin{figure}[htb]
\centering
\includegraphics[scale=1]{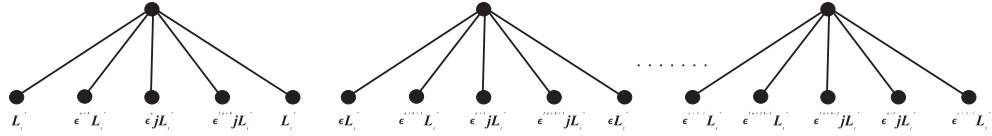}
\caption{The $n-4$ and $n-5$ level of the intersection poset.}
\label{fig:Fig7}
\end{figure}

\textit{(A)} $L_{1}^{\ast }\cap \epsilon ^{b+a}L_{1}^{\ast }\cap \epsilon
^{a}jL_{1}^{\ast }\cap \epsilon ^{2a+b}jL_{1}^{\ast }\cap L_{2}^{\ast }$ is
a linear space of dimension $n-5$

\textit{(B)} $L_{1}^{\ast }$ is a half-subspace of dimension $n-4$

\textit{(C)} $L_{2}^{\ast }$ is a linear space of dimension $n-4$

\noindent we conclude that%
\begin{eqnarray*}
H_{2}(W_{n}\backslash \cup \mathcal{A}(\mathcal{J},\alpha );\mathbb{Z})
&\cong &\mathbb{Z}^{a+b}\oplus \mathbb{Z}^{4(a+b)}\oplus \underset{d=0}{%
\overset{n-6}{\bigoplus }}\mathrm{Hom}\left( \underset{\dim
V=d}{\bigoplus
}H_{n-4}(\Delta (P_{<V})\ast \hat{V});\mathbb{Z}\right) \\
&\cong &\mathbb{Z}^{a+b}\oplus \mathbb{Z}^{4(a+b)}\oplus \underset{d=0}{%
\overset{n-6}{\bigoplus }}\mathrm{Hom}\left( \underset{\dim
V=d}{\bigoplus }\tilde{H}_{n-5-\dim V}(\Delta
(P_{<V}));\mathbb{Z}\right) .
\end{eqnarray*}

\begin{lemma}
$(\forall V\in P(\alpha ))\,\dim V\leq n-6\,\Longrightarrow \,\tilde{H}%
_{n-5-\dim V}(\Delta (P(\alpha )_{<V}))=0.$
\end{lemma}

\begin{proof}
For every element $W\in P(\alpha )_{<V}$ such that $\dim W=n-4$ there exists
a unique element $U_{W}\in P(\alpha )_{<V}$ with the property $\dim
U_{W}=n-5 $ and $W<U_{W}$. There is a monotone map $f:P(\alpha
)_{<V}\rightarrow P(\alpha )_{<V}-\{U~|~\dim U=n-4\}$ defined by
\begin{equation*}
\begin{array}{ccc}
U & \longmapsto & \left\{
\begin{array}{cc}
U, & \text{for }\dim U\leq n-5 \\
U_{W} & \text{for }\dim U=n-4%
\end{array}%
\right.%
\end{array}%
\end{equation*}%
which satisfies conditions of the Quillen fiber lemma. Thus, $f$ induces a
homotopy equivalence and so
\begin{equation*}
\tilde{H}_{n-5-\dim V}(\Delta (P(\alpha )_{<V}))=\tilde{H}_{n-5-\dim
V}(\Delta (P(\alpha )_{<V}-\{U~|~\dim U=n-4\}))=0.
\end{equation*}%
because $\dim \Delta (P(\alpha )_{<V}-\{U~|~\dim U=n-4\})<n-5-\dim V$.
\end{proof}

\begin{conclusion}
The previous lemma and equality impies that
\begin{equation}
H_{2}(W_{n}\backslash \cup \mathcal{A}(\mathcal{J},\alpha );\mathbb{Z})\cong
H_{n-4}(\cup \widehat{\mathcal{A}}(\mathcal{J},\alpha );\mathbb{Z})\cong
\mathbb{Z}^{a+b}\oplus \mathbb{Z}^{4(a+b)}.  \label{H2(complement)}
\end{equation}
\end{conclusion}

\subsection{Step 6}

\noindent We compute $\mathbb{Q}_{4n}$ coinvariants by working in the $%
\mathbb{Q}_{4n}$-module $H_{n-4}(\cup \widehat{\mathcal{A}}(\mathcal{J}%
,\alpha );\mathbb{Z})$ with the modified action. Since action respects
dimensional decomposition (\ref{hom decomposition - 1}) we analyze $n-4$ and
$n-5$ dimension cases separatelly.

\textit{(A)} Let $l\in H_{n-4}(\cup \widehat{\mathcal{A}}(\mathcal{J},\alpha
);\mathbb{Z})$ be the element which corresponds to the compactification of
the subspace $L_{2}^{\ast }$. Then $l,\epsilon l,..,\epsilon ^{a+b-1}l$ is a
base of the first factor in the equality (\ref{H2(complement)}). There are
to set identities which should be considered.

\noindent The identity $L_{2}^{\ast }=\epsilon ^{a+b}L_{2}^{\ast }$ produces
equality $l=\det (\epsilon ^{a+b})\epsilon ^{a+b}\cdot
l=(-1)^{(n+1)(a+b)}\epsilon ^{a+b}\cdot l$ in homology. Indeed, like in
examples \ref{Ex:Coinvariants1} \ and \ref{Ex:Coinvariants2}, we look at the
orthogonal complement of $L_{2}^{\ast }$ and see that $\epsilon ^{a+b}$ acts
on a basis as the even permutation $(3412)$. Then coinvariant calculation
with modified action implies a trivial relation,%
\begin{equation*}
l\sim \epsilon ^{a+b}\ast l=\det (\epsilon ^{a+b})\epsilon ^{a+b}\cdot
l=\det (\epsilon ^{a+b})^{2}l.
\end{equation*}%
The second identity $L_{2}^{\ast }=\epsilon ^{-b}jL_{2}^{\ast }$ produces
equality $l=(-1)\det (\epsilon ^{-b}j)\epsilon ^{-b}j\cdot l$ in homology.
Like we mentioned, this is the consequence of the fact that $\epsilon ^{-b}j$
acts on a base of $\left( L_{2}^{\ast }\right) ^{\bot }$ as the odd
permutation $(3214)$. Then,%
\begin{equation*}
l\sim \epsilon ^{-b}j\ast l=\det (\epsilon ^{-b}j)\epsilon ^{-b}j\cdot
l=-\det (\epsilon ^{-b}j)^{2}\cdot l=-l.
\end{equation*}

\noindent Thus, the first factor in the equality (\ref{H2(complement)})
reduces in the coinvariants to a single $\mathbb{Z}_{2}$.

\begin{figure}[htb]
\centering
\includegraphics[scale=1]{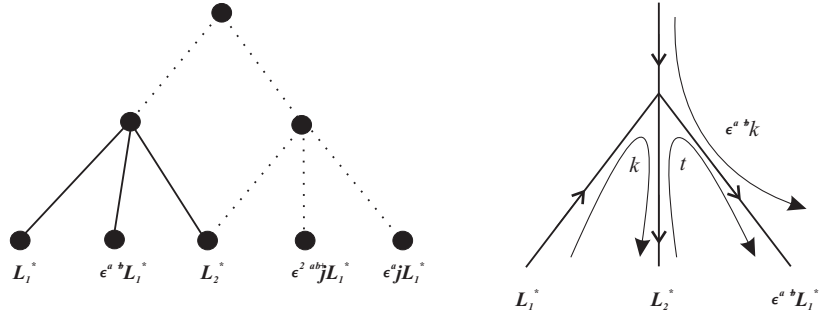}
\caption{Generating part of the Hasse diagram for coinvariants. }
\label{fig:Fig9}
\end{figure}

\textit{(B) }Since we want to find the coinvariants we can concentrate on
the "genterating" part of the Hasse diagram and its inner symmetries (figure %
\ref{fig:Fig9})


\noindent Let $k$ and $h$ be the elements of $H_{n-4}(\cup \widehat{\mathcal{%
A}}(\mathcal{J},\alpha );\mathbb{Z})$ which correspond to the
compactifications of the intersections $L_{1}^{\ast }\cap L_{2}^{\ast }$
and. $\epsilon ^{a+b}L_{1}^{\ast }\cap L_{1}^{\ast }$. Then the orbits of
these two elements form a base of the second factor in the equality (\ref%
{H2(complement)}). Let $t\in H_{n-4}(\cup \widehat{\mathcal{A}}(\mathcal{J}%
,\alpha );\mathbb{Z})$ also denotes the element which corresponds to the
compactification of the intersection $\epsilon ^{a+b}L_{1}^{\ast }\cap
L_{2}^{\ast }$ with the same generic point in $L_{2}^{\ast }$ as in the
representation of the $k$. Here we have to be very carefull, because $%
\epsilon ^{a+b}$ interchanges tha halfspaces of $L_{2}^{\ast }$
genterated by the hyperplane $H_{1}$. Moreover Thus, like the
figure \ref{fig:Fig9}
indicates%
\begin{equation}
\epsilon ^{a+b}\cdot k=\det (\epsilon ^{a+b})(l+t)  \label{rel1}
\end{equation}%
where $l\in H_{n-4}(\cup \widehat{\mathcal{A}}(\mathcal{J},\alpha );\mathbb{Z%
})$ (as in part (A)) is homology class corresponding to $\widehat{%
L_{2}^{\ast }}$.

\noindent On the other hand $k+t=h$. Since the sphere $\widehat{L_{1}^{\ast
}\cap \epsilon ^{a+b}L_{1}^{\ast }}\ast S^{0}$ is invariant of the element $%
\epsilon ^{a+b}$, there is a possible equality $h=\pm \epsilon ^{a+b}\cdot h$%
. Let us fix the basis $%
\{e_{1}+..+e_{a},~e_{a+1}+..+e_{a+b},~e_{a+b+1}+..+e_{2a+b},~e_{2a+b+1}+..+e_{n},~(a+b)(e_{a}-e_{2a+b}+e_{1}-e_{a+b+1})+e_{a+1}-e_{2a+b+1}+e_{n}-e_{a+b}\}
$ of the orthogonal complement of $L_{1}^{\ast }\cap \epsilon
^{a+b}L_{1}^{\ast }$. Then the matrix of $\epsilon ^{a+b}$ on this basis is%
\begin{equation*}
M=\left[
\begin{array}{lllll}
0 & 0 & 1 & 0 & 0 \\
0 & 0 & 0 & 1 & 0 \\
1 & 0 & 0 & 0 & 0 \\
0 & 1 & 0 & 0 & 0 \\
0 & 0 & 0 & 0 & -1%
\end{array}%
\right]
\end{equation*}%
and $\det M=-1$. This produces the relation $h=-\det (\epsilon
^{a+b})\epsilon ^{a+b}\cdot h$ in homology. Thus,%
\begin{equation}
h\sim \epsilon ^{a+b}\ast h=\det (\epsilon ^{a+b})\epsilon ^{a+b}\cdot
h=-\det (\epsilon ^{a+b})^{2}h=-h\text{.}  \label{rel2}
\end{equation}%
Since, $\epsilon ^{a+b}\cdot k-\det (\epsilon ^{a+b})t=\det (\epsilon
^{a+b})l$ and $h=k+t$, the relations \ref{rel1} and \ref{rel2} imply that%
\begin{equation*}
l\sim \det (\epsilon ^{a+b})l=\det (\epsilon ^{a+b})(\epsilon ^{a+b}\cdot
k-\det (\epsilon ^{a+b})t)=\epsilon ^{a+b}\ast k-t\sim k-t\text{ and }%
k+t\sim -k-t\text{.}
\end{equation*}%
Actually, we conclude that $4k=0$, or the second factor produces one $%
\mathbb{Z}_{4}$ .

\begin{conclusion}
$H_{2}(W_{n}\backslash \cup \mathcal{A}(\mathcal{J},\alpha );\mathbb{Z})_{%
\mathbb{Q}_{4n}}\cong H_{2}(W_{n}\backslash \cup \mathcal{A}(\mathcal{J}%
,\alpha );\mathbb{Z})_{\mathbb{Q}_{4n}}\cong \mathbb{Z}_{2}\oplus \mathbb{Z}%
_{4}$.
\end{conclusion}

\subsection{\label{Sec:Step7}Step 7}

\noindent Now we are ready to find the image of the obstruction cocycle $%
c^{\prime }=2\left\Vert v\right\Vert $ via the isomorphism $\vartheta
:H_{2}(W_{n}\backslash \cup \mathcal{A}(\mathcal{J},\alpha );\mathbb{Z}%
)\rightarrow \mathrm{Hom}\left( H_{n-4}(\cup \widehat{\mathcal{A}}(\mathcal{J%
},\alpha ),\mathbb{Z}),\mathbb{Z}\right) $. Actually we have to compute
(from (\ref{link-isomorphism})) the linking numbers of the sphere / boundary
of the simplex $S=\partial \lbrack u_{a},u_{a+1},u_{2a+b},u_{2a+b+1}]$
(which represents the broken point class $\left\Vert v\right\Vert $) and
spheres which represents generators of the $H_{n}(\cup \widehat{\mathcal{A}}(%
\mathcal{J},\alpha );\mathbb{Z})$. Since simplex $%
[u_{a},u_{a+1},u_{2a+b},u_{2a+b+1}]$ intersects only $I=L_{1}^{\ast }\cap
\epsilon ^{b+a}L_{1}^{\ast }\cap \epsilon ^{a}jL_{1}^{\ast }\cap \epsilon
^{2a+b}jL_{1}^{\ast }\cap L_{2}^{\ast }$, the homology elements which may
have non-zero values are those from the local picture $H_{n}(\widehat{I}\ast
\lbrack 5];\mathbb{Z})$. Since the coinvariant class of $l$ lives in $%
\mathbb{Z}_{2}$ and our obstruction cocycle is $2\cdot (something)$ there is
no need to compute the image of $l$ either.

\noindent By moving simplex $[u_{a},u_{a+1},u_{2a+b},u_{2a+b+1}]$ and
analyzing its intersection with subspaces $L_{1}^{\ast }$, $\epsilon
^{b+a}L_{1}^{\ast }$, $\epsilon ^{a}jL_{1}^{\ast }$, $\epsilon
^{2a+b}jL_{1}^{\ast }$, $L_{2}^{\ast }$ we decompose the broken point class $%
\left\Vert v\right\Vert $ into a sum of ordinary point classes.

\begin{theorem}
\label{Lemma:Link}Let $v_{1}\in L_{1}^{\ast }$, $v_{2}\in \epsilon
^{b+a}L_{1}^{\ast }$, $v_{3}\in \epsilon ^{a}jL_{1}^{\ast }$, $v_{4}\in
\epsilon ^{2a+b}jL_{1}^{\ast }$ and $y_{1}$, $y_{2}\in L_{2}^{\ast }$ are
arbitrary elements with the property that thay do not belong to any other
element of the arrangement $\mathcal{A}(\mathcal{J},\alpha )$. The request
for $y_{1}$ and $y_{2}$ is that thay are in different connecting components
of $L_{2}^{\ast }\backslash I$. Then our broken point class $\left\Vert
v\right\Vert $ can be presented as the sum of the point classes%
\begin{equation}
\left\Vert v\right\Vert =\tau _{1}\left\Vert v_{1}\right\Vert +\tau
_{2}\left\Vert v_{3}\right\Vert +\tau _{3}\left\Vert y_{1}\right\Vert =\mu
_{1}\left\Vert v_{2}\right\Vert +\mu _{2}\left\Vert v_{4}\right\Vert +\mu
_{3}\left\Vert y_{2}\right\Vert  \label{Broken-point}
\end{equation}%
where $\tau _{i}$ and $\mu _{i}$ are appropriate signs \{+1,-1\}.
\end{theorem}

\begin{proof}
The idea of the proof is very simple. \textit{Move our filled sphere /
simplex a little and look where it hits our arrangement and we are done.}
The reason for this to work is in the definition of the linking -
intersection number and the deffinition of (broken) point classes. Let us
recall that the linking /intersection number is "very" invariant under small
movement of the sphere - disc. Thus if you want to check whether you
calculated the intersection number properly you move your sphere - disc and
see what happens. Since we are in linear / convex situation, besides the
linking - intersection numbers (if correctly defined) are +1,-1 or 0, it is
enought to consider small translatorions of a sphere - disc.

\noindent Thus, let us move our simplex $\sigma $ by the generic
"small" vector $s=\sum_{i=1}^{n}\xi _{i}e_{i}$. Because of
complementary dimension affine space
$\mathrm{span}\{u_{a+b},u_{a+b+1},u_{n},u_{1}\}+s$ must hit all
linear spans of $L_{1}^{\ast }$, $\epsilon ^{a+b}L_{1}^{\ast }$,
$\epsilon ^{a}jL_{1}^{\ast }$, $\epsilon ^{2a+b}jL_{1}^{\ast }$
and $L_{2}^{\ast }$. Let us denote these intersection points by
\begin{eqnarray*}
v_{1} &\in &\mathrm{span}\left( L_{1}^{\ast }\right) \cap (\sigma +s)\text{,
}v_{2}\in \mathrm{span}\left( \epsilon ^{a+b}L_{1}^{\ast }\right) \cap
(\sigma +s)\text{, }v_{3}\in \mathrm{span}\left( \epsilon ^{a}jL_{1}^{\ast
}\right) \cap (\sigma +s)\text{, } \\
v_{4} &\in &\mathrm{span}\left( \epsilon ^{2a+b}jL_{1}^{\ast }\right) \cap
(\sigma +s)\text{, }w\in \mathrm{span}\left( L_{2}^{\ast }\right) \cap
(\sigma +s).
\end{eqnarray*}%
If the translation is small enough these points will remain in the interion
of the simplex $\sigma +s$. Now the tiresome, but necessary part. The points
are given by%
\begin{eqnarray*}
{\small v}_{1} &{\small =}&\tfrac{\alpha _{1}}{n}{\small e}_{a}{\small +}%
\tfrac{\beta _{1}}{n}{\small e}_{a+1}{\small +}\tfrac{\gamma _{1}}{n}{\small %
e}_{2a+b}{\small +}\tfrac{\delta _{1}}{n}{\small e}_{2a+b+1}{\small -}\tfrac{%
1}{n}\sum_{i=1}^{n}{\small e}_{i}{\small +}\sum_{i=1}^{n}{\small \xi }_{i}%
{\small e}_{i}, \\
{\small \alpha }_{1} &{\small =}&{\small a-n}\sum_{1}^{a}{\small \xi }_{i}%
{\small ,~\gamma _{1}=a+b-n\sum_{a+1}^{2a+b}\xi _{i}-\beta
_{1},~\delta
_{1}=b-n}\sum_{2a+b+1}^{n}{\small \xi }_{i}, \\
{\small \beta }_{1} &{\small =}&\tfrac{{\small -n}}{a+b+1}\left(
\sum_{2a+b+2}^{n}{\small \xi }_{i}+{\small \xi }_{n}-{\small \xi }_{a+b}+%
{\small \xi }_{a+1}{\small -}\left( {\small a+b}\right) \left(
\sum_{2}^{a-1}{\small \xi }_{i}{\small -}\sum_{a+1}^{2a+b-1}{\small \xi }%
_{i}+{\small \xi }_{a+b+1}\right) \right) {\small +b}
\end{eqnarray*}%
\begin{eqnarray*}
{\small v}_{2} &{\small =}&\tfrac{\alpha _{2}}{n}{\small e}_{a}{\small +}%
\tfrac{\beta _{2}}{n}{\small e}_{a+1}{\small +}\tfrac{\gamma _{2}}{n}{\small %
e}_{2a+b}{\small +}\tfrac{\delta _{2}}{n}{\small e}_{2a+b+1}{\small -}\tfrac{%
1}{n}\sum_{i=1}^{n}{\small e}_{i}{\small +}\sum_{i=1}^{n}{\small \xi }_{i}%
{\small e}_{i}, \\
{\small \alpha }_{2} &{\small =}&{\small
a+b-n}\sum_{2a+b+1}^{n}{\small \xi
}_{i}{\small -n}\sum_{i=1}^{a}{\small \xi }_{i}-{\small \delta }_{2},~%
{\small \beta }_{2}{\small =b-n}\sum_{a+1}^{a+b}{\small \xi }_{i}{\small %
,~\gamma }_{2}{\small =a-n}\sum_{a+b+1}^{2a+b}{\small \xi }_{i}, \\
{\small \delta }_{2} &{\small =}&\tfrac{{\small -n}}{a+b+1}\left(
\sum_{a+2}^{a+b}{\small \xi }_{i}{\small -\xi }_{n}{\small +\xi }_{a+b}%
{\small +\xi }_{2a+b+1}+\left( {\small a+b}\right) \left( \sum_{2a+b+1}^{n}%
{\small \xi }_{i}{\small +}\sum_{2}^{a-1}{\small \xi }_{i}{\small -}%
\sum_{a+b+2}^{2a+b-1}{\small \xi }_{i}\right) \right) {\small +b;}
\end{eqnarray*}

\begin{eqnarray*}
{\small v}_{3} &{\small =}&\tfrac{\alpha _{3}}{n}{\small e}_{a}{\small +}%
\tfrac{\beta _{3}}{n}{\small e}_{a+1}{\small +}\tfrac{\gamma _{3}}{n}{\small %
e}_{2a+b}{\small +}\tfrac{\delta _{3}}{n}{\small e}_{2a+b+1}{\small -}\tfrac{%
1}{n}\sum_{i=1}^{n}{\small e}_{i}{\small +}\sum_{i=1}^{n}{\small \xi }_{i}%
{\small e}_{i}, \\
{\small \alpha }_{3} &{\small =}&{\small a-n}\sum_{1}^{a}{\small \xi }_{i},~%
{\small \beta }_{3}{\small =b-n}\sum_{a+1}^{a+b}{\small \xi }_{i}{\small %
,~\gamma }_{3}{\small =a+b-n}\sum_{a+b+1}^{n}{\small \xi }_{i}{\small %
-\delta }_{3} \\
{\small \delta }_{3} &{\small =}&{\small \tfrac{-n}{a+b-1}\left(
(a+b)\left( \sum_{a+b+2}^{n}\xi _{i}-\sum_{2}^{a-1}\xi _{i}-\xi
_{2a+b}\right) -\sum_{a+2}^{a+b}\xi _{i}-\xi _{2a+b+1}-\xi
_{a+b}+\xi _{n}\right) +b;}
\end{eqnarray*}

\begin{eqnarray*}
{\small v}_{4} &{\small =}&\tfrac{\alpha _{4}}{n}{\small e}_{a}{\small +}%
\tfrac{\beta _{4}}{n}{\small e}_{a+1}{\small +}\tfrac{\gamma _{4}}{n}{\small %
e}_{2a+b}{\small +}\tfrac{\delta _{4}}{n}{\small e}_{2a+b+1}{\small -}\tfrac{%
1}{n}\sum_{i=1}^{n}{\small e}_{i}{\small +}\sum_{i=1}^{n}{\small \xi }_{i}%
{\small e}_{i}, \\
{\small \alpha }_{4} &{\small =}&{\small a+b-n}\sum_{1}^{a+b}{\small \xi }%
_{i}{\small -\beta }_{4}{\small ,}~{\small \delta
_{4}=b-n\sum_{2a+b+1}^{n}\xi _{i},}~{\small \gamma _{4}=a-n}%
\sum_{a+b+1}^{2a+b}{\small \xi }_{i} \\
{\small \beta }_{4} &{\small =}&\tfrac{n}{a+b-1}\left( {\small (a+b)}\left(
\sum_{a+b+1}^{2a+b}{\small \xi }_{i}{\small -}\sum_{2}^{a+b}{\small \xi }%
_{i}{\small +\xi }_{a}\right) -\sum_{2a+b+1}^{n-1}{\small \xi }_{i}{\small %
-\xi }_{2a+b+1}{\small +\xi }_{a+1}{\small -\xi }_{a+b}\right) {\small +b;}
\end{eqnarray*}

\begin{eqnarray*}
{\small w}_{4} &{\small =}&\tfrac{\alpha _{5}}{n}{\small e}_{a}{\small +}%
\tfrac{\beta _{5}}{n}{\small e}_{a+1}{\small +}\tfrac{\gamma _{5}}{n}{\small %
e}_{2a+b}{\small +}\tfrac{\delta _{5}}{n}{\small e}_{2a+b+1}{\small -}\tfrac{%
1}{n}\sum_{i=1}^{n}{\small e}_{i}{\small +}\sum_{i=1}^{n}{\small \xi }_{i}%
{\small e}_{i}, \\
{\small \alpha }_{5} &=&{\small a-n\sum_{1}^{a}\xi _{i},}~{\small \beta }%
_{4}{\small =b-n\sum_{a+1}^{a+b}\varepsilon _{i},~\gamma _{5}=a-n}%
\sum_{a+b+1}^{2a+b}{\small \xi }_{i},~{\small \delta
_{5}=b-n\sum_{2a+b+1}^{n}\xi _{i}}\text{.}
\end{eqnarray*}

\noindent Now we will prove the following equivalences%
\begin{eqnarray*}
&&\left( v_{1}\in K^{+}\Leftrightarrow v_{2}\in K^{+}\right) ~\vee ~\left(
v_{1}\in K^{-}\Leftrightarrow v_{2}\in K^{-}\right)  \\
&&\left( v_{3}\in \epsilon ^{a}jK^{+}\Leftrightarrow v_{4}\in \epsilon
^{a}jK^{+}\right) ~\vee ~\left( v_{3}\in \epsilon ^{a}jK^{-}\Leftrightarrow
v_{4}\in \epsilon ^{a}jK^{-}\right)
\end{eqnarray*}%
which have the following consequence
\begin{equation}
\left( v_{1}\in L_{1}^{\ast }\Leftrightarrow v_{2}\notin \epsilon
^{a+b}L_{1}^{\ast }\right) \text{ and }\left( v_{3}\in \epsilon
^{a}jL_{1}^{\ast }\Leftrightarrow v_{4}\notin \epsilon ^{2a+b}jL_{1}^{\ast
}\right)   \label{equvalencije}
\end{equation}%
and therefore the equality (\ref{Broken-point}).stands.

\noindent In order to prove equivalences (\ref{equvalencije}) we evaluate $%
v_{1}$, $v_{2}$, $v_{3}$ and $v_{4}$\ on linear forms $x_{1}+...+x_{a+b}$, $%
x_{a+b+1}+...+x_{n}$, $x_{2a+b+1}+..+x_{n}+x_{1}+..+x_{a}$ and $%
x_{a+1}+..+x_{2a+b}$, respectively. It can be calculated that
\begin{eqnarray*}
&&%
\begin{array}{ll}
{\small x}_{1}{\small +..+x}_{a+b}{\small ~|}_{v_{1}} & {\small =x}_{a+1}%
{\small +..+x}_{a+b}{\small ~|}_{v_{1}}{\small =}\tfrac{\beta _{1}}{n}%
{\small -}\tfrac{b}{n}{\small +}\sum_{a+1}^{a+b}{\small \xi }_{i} \\
& {\small =}\tfrac{a+b}{a+b+1}\sum_{2}^{a-1}{\small \xi }_{i}{\small +}%
\tfrac{1}{a+b+1}\sum_{a+2}^{a+b-1}{\small \xi }_{i}{\small +}\tfrac{2}{a+b+1%
}{\small \xi }_{a+b} \\
& {\small -}\tfrac{a+b}{a+b+1}\sum_{a+b+2}^{2a+b-1}{\small \xi }_{i}{\small %
-}\tfrac{1}{a+b+1}\sum_{2a+b+2}^{n-1}{\small \xi }_{i}{\small -}\tfrac{2}{%
a+b+1}{\small \xi }_{n}%
\end{array}
\\
&&%
\begin{array}{ll}
{\small x}_{a+b+1}{\small +..+x}_{n}{\small ~|}_{v_{2}} & {\small =x}%
_{2a+b+1}{\small +...+x}_{n}{\small ~|}_{v_{2}}{\small =}\tfrac{\delta _{2}}{%
n}{\small -}\tfrac{b}{n}{\small +}\sum_{2a+b+1}^{n}{\small \xi }_{i} \\
& ={\small -}\tfrac{a+b}{a+b+1}\sum_{i=2}^{a-1}{\small \xi }_{i}{\small -}%
\tfrac{1}{a+b+1}\sum_{a+2}^{a+b-1}{\small \xi }_{i}{\small -}\tfrac{2}{a+b+1%
}{\small \xi }_{a+b} \\
& {\small +}\tfrac{a+b}{a+b+1}\sum_{i=a+b+2}^{2a+b-1}{\small \xi }_{i}%
{\small +}\tfrac{1}{a+b+1}\sum_{2a+b+2}^{n-1}{\small \xi }_{i}{\small +}%
\tfrac{2}{a+b+1}{\small \xi }_{n}%
\end{array}
\\
&&%
\begin{array}{ll}
{\small x}_{2a+b+1}{\small +..+x}_{n}{\small +x}_{1}{\small +..+x}_{a}%
{\small ~|}_{v_{3}} & {\small =x}_{2a+b+1}{\small +...+x}_{n}{\small ~|}%
_{v_{2}}{\small =}\tfrac{\delta _{3}}{n}{\small -}\tfrac{b}{n}{\small +}%
\sum_{2a+b+1}^{n}{\small \xi }_{i} \\
& =\tfrac{a+b}{a+b-1}\sum_{2}^{a-1}{\small \xi }_{i}{\small +}\tfrac{1}{%
a+b-1}\sum_{a+2}^{a+b-1}{\small \xi }_{i}{\small +}\tfrac{2}{a+b-1}{\small %
\xi }_{a+b} \\
& {\small -}\text{{\small $\tfrac{a+b}{a+b-1}\sum_{a+b+2}^{2a+b-1}\xi _{i}-%
\tfrac{1}{a+b-1}\sum_{2a+b+2}^{n-1}\xi _{i}-\tfrac{2}{a+b-1}\xi _{n}$}}%
\end{array}
\\
&&%
\begin{array}{ll}
x_{a+1}+..+x_{2a+b}{\small ~|}_{v_{4}} & {\small =x}_{a+1}{\small +...+x}%
_{a+b}{\small ~|}_{v_{4}}{\small =}\tfrac{\beta _{4}}{n}{\small -}\tfrac{b}{n%
}{\small +\sum_{a+1}^{a+b}\xi _{i}} \\
& =-\tfrac{a+b}{a+b-1}\sum_{2}^{a-1}{\small \xi }_{i}{\small -}\tfrac{1}{%
a+b-1}\sum_{a+2}^{a+b-1}{\small \xi }_{i}{\small -}\tfrac{2}{a+b-1}{\small %
\xi }_{a+b} \\
& {\small +}\text{{\small $\tfrac{a+b}{a+b-1}\sum_{a+b+2}^{2a+b-1}\xi _{i}+%
\tfrac{1}{a+b-1}\sum_{2a+b+2}^{n-1}\xi _{i}+\tfrac{2}{a+b-1}\xi _{n}$}}%
\end{array}%
\end{eqnarray*}

\noindent Thus, we have proved
\begin{eqnarray*}
\left( {\small a+b+1}\right) \left( {\small x}_{1}{\small +..+x}_{a+b}%
{\small ~|}_{v_{1}}\right) &=&-\left( {\small a+b+1}\right) \left( {\small x}%
_{a+b+1}{\small +...+x}_{n}{\small ~|}_{v_{2}}\right) \\
&=&\left( {\small a+b-1}\right) \left( {\small x}_{2a+b+1}{\small +..+x}_{n}%
{\small +x}_{1}{\small +..+x}_{a}{\small ~|}_{v_{3}}\right) \\
&=&-\left( {\small a+b+1}\right) \left( x_{a+1}+..+x_{2a+b}{\small ~|}%
_{v_{4}}\right)
\end{eqnarray*}%
which implies the equivalence (\ref{equvalencije}).

\noindent It looks like we forgot about the space $L_{2}^{\ast }$ and point $%
w$. The reason is that when ever we move the simplex $\sigma $ it will hit
the space $L_{2}^{\ast }$. In respect of the direction we move $\sigma $ the
point $w$ will be in one or in another connecting component of $L_{2}^{\ast
}-I$. Thus, the broken point class $\left\Vert w\right\Vert $ appears in
both equalities. Also, observe that we do not worry about the signs $\tau
_{i}$ and $\mu _{i}$, they will not play any role in our calculation.


\begin{figure}[htb]
\centering
\includegraphics[scale=0.8]{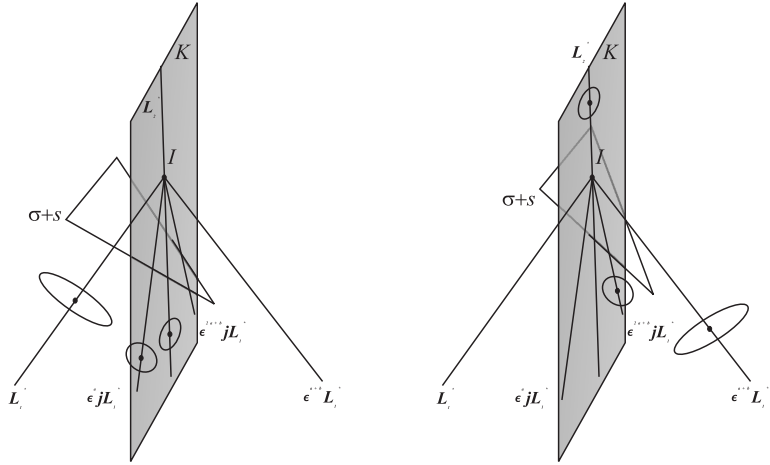}
\caption{Generating part of the Hasse diagram for coinvariants. }
\label{fig:Fig8}
\end{figure}

\end{proof}

\begin{conclusion}
Let us observe that $\left\Vert v_{1}\right\Vert =\pm \epsilon
^{a}j\left\Vert v_{3}\right\Vert $ and $\left\Vert v_{2}\right\Vert =\pm
\epsilon ^{a}j\left\Vert v_{4}\right\Vert $. Thus, since the class of our
obstruction cocycle%
\begin{equation*}
c^{\prime }=2\left\Vert v\right\Vert =2\left( \tau _{1}\left\Vert
v_{1}\right\Vert +\tau _{2}\left\Vert v_{3}\right\Vert +\tau _{3}\left\Vert
y_{1}\right\Vert \right) =2\left( \mu _{1}\left\Vert v_{2}\right\Vert +\mu
_{2}\left\Vert v_{4}\right\Vert +\mu _{3}\left\Vert y_{2}\right\Vert \right)
\end{equation*}%
lives in coinvariants $H_{2}(W_{n}\backslash \cup \mathcal{A}(\mathcal{J}%
,\alpha );\mathbb{Z})_{\mathbb{Q}_{4n}}\cong \mathbb{Z}_{2}\oplus \mathbb{Z}%
_{4}$ we can shift our attention to the cohomological cocycle%
\begin{equation*}
c^{\prime \prime }=2\tau _{3}\left\Vert y_{1}\right\Vert .
\end{equation*}
\end{conclusion}

\subsection{\label{sec:Step8}Step 8}

\noindent Before proving that the cohomology class of the obstruction
cocycle $c^{\prime \prime }=2\tau _{3}\left\Vert y_{1}\right\Vert $ is not
zero, let us recall the nature of the Poincar\'{e} duality isomorphism $%
\vartheta :H_{m-1}(\mathbb{R}^{n+m}\backslash \cup \mathcal{A};\mathbb{Z}%
)\rightarrow \mathrm{Hom}(H_{n}(\cup \widehat{\mathcal{A}};\mathbb{Z}),%
\mathbb{Z})$. As we have seen, for $t\in H_{n}(\cup \widehat{\mathcal{A}};%
\mathbb{Z});\mathbb{Z})$ represented by a submanifold $T$ in $\cup \widehat{%
\mathcal{A}}$, the image $\vartheta (\left\Vert x\right\Vert )$ of $t$ is
given by the linking number%
\begin{equation*}
\vartheta (\left\Vert x\right\Vert )(t)=\mathrm{link}(S,T)\text{..}
\end{equation*}

\begin{theorem}
As we have already denoted, let $l,k,h\in H_{n-4}(\cup \widehat{\mathcal{A}}(%
\mathcal{J},\alpha );\mathbb{Z})$ be elements corresponding to spheres $%
\widehat{L_{2}^{\ast }}$, $\widehat{L_{1}^{\ast }\cap L_{2}^{\ast }}\ast
S^{\ast }$, $\widehat{\epsilon ^{a+b}L_{1}^{\ast }\cap L_{1}^{\ast }}\ast
S^{\ast }$. Then there exists a basis of $H_{n-4}(\cup \widehat{\mathcal{A}}(%
\mathcal{J},\alpha );\mathbb{Z})$ which contains $l$, $k$, and $h$ such that
\begin{equation}
\vartheta (\left\Vert y_{1}\right\Vert )(l)=\pm 1,\vartheta (\left\Vert
y_{1}\right\Vert )(k)=\pm 1,\vartheta (\left\Vert y_{1}\right\Vert )(h)=0
\label{linkovanje}
\end{equation}%
and for every other base element $d\in $ $H_{n-4}(\cup \widehat{\mathcal{A}}(%
\mathcal{J},\alpha );\mathbb{Z})$, $\vartheta (\left\Vert y_{1}\right\Vert
)(d)=0$. Thus,

(1) $\vartheta (\left\Vert y_{1}\right\Vert )=\mu _{1}l+\mu _{2}k\in
H_{n-4}(\cup \widehat{\mathcal{A}}(\mathcal{J},\alpha );\mathbb{Z})$ where $%
\mu _{i}\in \{+1,-1\}$;

(2) $c^{\prime \prime }=2k\neq 0$ in coinvariants $H_{2}(W_{n}\backslash
\cup \mathcal{A}(\mathcal{J},\alpha );\mathbb{Z})_{\mathbb{Q}_{4n}}\cong
H_{n}(\cup \widehat{\mathcal{A}};\mathbb{Z})_{\mathbb{Q}_{4n}}$.

\noindent We abuse notation in (2) and (3) using the fact that $H_{n-4}(\cup
\widehat{\mathcal{A}}(\mathcal{J},\alpha );\mathbb{Z})$ is free ($\mathrm{Hom%
}$ can be deleted) and omitting the notation of the cohomology class to
simplify notation.
\end{theorem}

\begin{proof}
The proof can be read from the figure \ref{fig:Fig8}. For example,
the required basis of $H_{n-4}(\cup
\widehat{\mathcal{A}}(\mathcal{J},\alpha );\mathbb{Z}) $ is
composed of

(A) the orbit of $l$,

(B) $k$, $h$, $\epsilon ^{a+b}h$ and element $f$ corresponding to the sphere
$\widehat{\epsilon ^{a+b}L_{1}^{\ast }\cap \epsilon ^{a}jL_{1}^{\ast }},$

(C) anything you want outside the $I$ "neighbourhood", or precisely $%
\epsilon k$, $\epsilon h$, $\epsilon ^{a+b+1}h$, $\epsilon f$, ..., $%
\epsilon ^{a+b-1}k$, $\epsilon ^{a+b-1}h$, $\epsilon ^{n-1}h$, $\epsilon
^{a+b-1}f$.

\noindent This basis satisfies the property (\ref{linkovanje}). Thus (1) and
(2) follow directly with the knowing that $4k=0$ in coinvariants.
\end{proof}

\begin{conclusion}
Since we proved that the cohomology class of the obstruction cocycle is not
zero we are ready to sum up (remember Propositions \ref%
{prop:PrelazakNaNovuGrupu} and \ref{prop:VezaProblem-Ekvivarijantan}):
\begin{equation*}
\begin{tabular}{lll}
$\text{There is no }\mathbb{Q}_{4n}\text{-map }S^{3}\rightarrow \cup
\mathcal{A}(\mathcal{J}\text{,}\alpha )\text{ }$ & $\Longrightarrow $ &
There is no $\mathbb{Q}_{4n}$-map $S^{3}\rightarrow \cup \mathcal{A}$ \\
& $\Longleftrightarrow $ & There is no $\mathbb{D}_{2n}$-map $V_{2}(\mathbb{R%
}^{3})\rightarrow \cup \mathcal{A}$ \\
& $\Longrightarrow $ & There is an $(\frac{a}{n},\frac{a+b}{n},\frac{b}{n})$%
-partition.%
\end{tabular}%
\text{\ }
\end{equation*}%
Once more, we have just proved that for every $a,b>1$ and every two measures
$\mu $ and $\nu $ on $S^{2}$, there exists an $(\frac{a}{n},\frac{a+b}{n},%
\frac{b}{n})$-partition $(x;l_{1},l_{2},l_{3})$ of measures $\mu $ and $\nu $%
.
\end{conclusion}

\subsection{The limit argument}

\noindent The reasons for assuming nice properties for our measures will
finally pop up. Let $\emph{S}\subseteq \mathbb{R}_{>0}^{3}$ be the space of
all $3$-fan partitions, i.e. the space of all triples $\alpha =(a,b,c)$, $%
a+b+c=1$ such that there exists an $\alpha $-partition of measures $\mu $
and $\nu $ by a $3$-fan. Since our measures $\mu $ and $\nu $ are proper
Borel probability measures, the space $\emph{S}$ is a closed subset of $%
\mathbb{R}_{>0}^{3}$. All we did so far is proving the following inclusion
\begin{equation*}
\{(\tfrac{a}{n},\tfrac{a+b}{n},\tfrac{b}{n})\in \mathbb{Q}%
_{>0}^{3}~|~2a+2b=n,~a,b\in \mathbb{Z}\}\subseteq \emph{S.}
\end{equation*}%
The fact that $\emph{S}$ is closed implies the inclusion%
\begin{equation*}
\{(a,b,c)\in \mathbb{R}_{>0}^{3}~|~a+b+c=1\}=\mathrm{cl}\{(\tfrac{a}{n},%
\tfrac{a+b}{n},\tfrac{b}{n})~|~2a+2b=n,~a,b\in \mathbb{Z}\}\subseteq \emph{S}%
\text{.}
\end{equation*}%
and therefore we are finally done.

\end{document}